\documentstyle[amssymb,amsfonts]{amsart}

\def\R{{\hbox{\bf R}}}
\def\C{{\hbox{\bf C}}}

\def\N{{\hbox{\bf N}}}
\def\I{{\hbox{\bf I}}}

\font \roman = cmr10 at 10 true pt

\def\be#1{\begin{equation}\label{#1}}
\def\bas{\begin{align*}}
\def\eas{\end{align*}}
\def\bi{\begin{itemize}}
\def\ei{\end{itemize}}

\def\dist{{\hbox{\roman dist}}}

\def\eps{\varepsilon}
\newenvironment{proof}{\noindent {\bf Proof} }{\endprf\par}
\def \endprf{\hfill  {\vrule height6pt width6pt depth0pt}\medskip}
\def\emph#1{{\it #1}}
\def\textbf#1{{\bf #1}}

\theoremstyle{plain}
  \newtheorem{theorem}[subsection]{Theorem}
  \newtheorem{conjecture}[subsection]{Conjecture}

  \newtheorem{proposition}[subsection]{Proposition}
  
  \newtheorem{lemma}[subsection]{Lemma}
  \newtheorem{corollary}[subsection]{Corollary}

\theoremstyle{remark}

\theoremstyle{definition}
  \newtheorem{definition}[subsection]{Definition}

\begin{document}

\title{$L^p$ improving bounds for averages along curves}

\author{Terence Tao}
\address{Department of Mathematics, UCLA, Los Angeles CA 90095-1555}
\email{tao@@math.ucla.edu}

\author{James Wright}
\address{School of Mathematics, University of
Edinburgh, JCMB, King's Buildings, Mayfield Road, Edinburgh EH9 3JZ,
Scotland}
\email{wright@@maths.ed.ac.uk}



\begin{abstract}  We establish local $(L^p,L^q)$ mapping properties for averages on curves.  The exponents are sharp except for endpoints.
\end{abstract}

\maketitle

\section{Introduction}

Let $n \geq 2$, and let $M_1$ and $M_2$ be two smooth $n-1$-dimensional manifolds\footnote{Our conventions may appear somewhat strange, but this choice of notation will be convenient to work with later on.}, each containing a preferred origin $0_{M_1}$ and $0_{M_2}$.  We shall abuse notation and write $0$ for both $0_{M_1}$ and $0_{M_2}$.  For the purposes of integration we shall place a smooth Riemannian metric on $M_1$ and $M_2$, although the exact choice of this metric will not be relevant.  All our considerations shall be local to the origin 0.  

We are interested in the local $L^p$ improving properties of averaging operators on curves.  Before we give the rigorous description of these operators, let us first give an informal discussion.

Informally, we assume that we have a smooth assignment $x_2 \mapsto \gamma_{x_2}$ taking points in $M_2$ to curves in $M_1$, with a corresponding dual assignment $x_1 \to \gamma^*_{x_1}$ taking points in $M_1$ to curves in $M_2$, such that
$$ x_1 \in \gamma_{x_2} \iff x_2 \in \gamma^*_{x_1}.$$
We then form the operator $R$ taking functions\footnote{All our functions in this paper will be real-valued.} on $M_1$ to functions on $M_2$, defined by
\be{r-intuitive}
Rf(x_2) := \int_{\gamma_{x_2}} f(x_1) a(x_1,x_2)
\end{equation}
where the amplitude function $a(x_1,x_2)$ is just a smooth cutoff to a neighborhood of $(x_1,x_2) = (0,0)$, and integration on $\gamma_{x_2}$ is respect to the induced Riemannian metric from $M_2$.  This operator has an adjoint
$$ R^* g(x_1) := \int_{\gamma^*_{x_1}} g(x_2) a^*(x_1,x_2)$$
where the amplitude function $a^*$ is another smooth cutoff to $(x_1,x_2) = (0,0)$.
To keep things from being vacuous we shall assume that $0 \in \gamma_0$ (or equivalently that $0 \in \gamma^*_0$).

A basic example of such an operator is the two-dimensional \emph{Radon transform}, in which $n=3$, $M_1$ is the plane $\R^2$, and $M_2$ is the space of lines in $\R^2$ (endowed with some reasonable Riemannian metric, and with some line through the origin designated as the 0 line), with $\gamma$ given by the tautological map $\gamma_l := l$.  The dual curve $\gamma_x$ then consists of those lines in $M_2$ which contain $x$.  The Radon transform is then the familiar operator
$$ Rf(l) = \int_l f.$$

Another example of such an operator is that of convolution with a fixed curve $\gamma^*_0 \subset \R^{n-1}$, for instance $\gamma^*_0 := \{ (t,t^2, \ldots, t^{n-1}): t \in \R \}$.  We then define $M_1 := M_2 := \R^{n-1}$ and set $\gamma_{x_2} := x_2 - \gamma^*_0$, so the dual curve is $\gamma^*_{x_1} = x_1 + \gamma^*_0$.  We can then consider operators $R$ of the form 
$$ Rf = f * d\sigma$$
where $d\sigma$ is a fixed smooth compactly supported measure on $\gamma^*_0$.

This class of operators was first introduced in \cite{gelfand-2}, see also
\cite{guillemin-1}, \cite{hormander}, \cite{helgason}.
We are interested in determining the exponents $1 \leq p, q \leq \infty$ such that $R$ maps $L^p(M_1)$ to $L^q(M_2)$ locally in a neighborhood of the origin, perhaps shrinking the support of the cutoff $a$ if necessary.  As will become clear, such a mapping property is easy to show if $q \leq p$; the interesting case is when $q > p$.  Observe that such a mapping property is essentially independent of the metric structure on $M_1$, $M_2$ or on the cutoff $a$ (assuming the support is sufficiently small), and is indeed invariant under diffeomorphisms of $M_1$ and $M_2$.  Thus we expect the set of exponents $(p,q)$ with this mapping property to depend only on diffeomorphism-invariant properties of the family of curves $x_2 \mapsto \gamma_{x_2}$.

Having given an informal description of our problem of interest, we now give a formal setup (based on that in \cite{gelfand-2}, \cite{guillemin-1}; a prototype of this ``double fibration'' formulation appeared earlier in \cite{helgason}) which will be more convenient to work with.  In addition to the $n-1$-dimensional manifolds $M_1$ and $M_2$ given earlier, we shall also work with an $n$-dimensional smooth Riemannian manifold $\Sigma$ with another preferred origin $0 = 0_\Sigma$.  We will also assume we have smooth maps $\pi_1: \Sigma \to M_1$ and $\pi_2: \Sigma \to M_2$ which map $0$ to $0$, and are submersions (i.e. the derivative has maximal rank $n-1$) in a neighborhood of $0$.

We then define the transform $R$ by duality as
\be{duality}
\begin{split}
\int_{M_2} Rf(x_2) g(x_2)\ dx_2 &= \int_{M_1} f(x_1) R^* g(x_1)\ dx_1\\
&:= \int_\Sigma f(\pi_1(x)) g(\pi_2(x)) a(x)\ dx
\end{split}
\end{equation}
where $dx$ is the measure on $\Sigma$ induced by the Riemannian metric, and $a(x)$ is a smooth cutoff to a neighborhood of $0$.  (We shall always use $x$ to denote elements of $\Sigma$, and $x_1$, $x_2$ to denote elements of $M_1$, $M_2$ respectively).

To connect the above formalism with the intuitive discussion given earlier, we assume that the assignment $x_2 \mapsto \gamma_{x_2}$ and the amplitude function $a(x_1,x_2)$ is given, and then define $\Sigma$ to be the $n$-dimensional manifold
$$ \Sigma := \{ (x_1,x_2) \in M_1 \times M_2: x_1 \in \gamma_{x_2} \}
= \{ (x_1,x_2) \in M_1 \times M_2: x_2 \in \gamma^*_{x_1} \}$$
with $0_\Sigma := (0_{M_1}, 0_{M_2})$.  We then define $\pi_1$ and $\pi_2$ 
to be just the co-ordinate projections $\pi_1(x_1,x_2) := x_1$, $\pi_2(x_1,x_2) := x_2$.  The reader may easily verify by the change-of-variables formula that the transform given by \eqref{r-intuitive} then obeys the formula \eqref{duality}, if the support of the cutoffs are sufficiently small.  In particular the examples of the Radon transform and convolutions with curves given previously can be put into the above framework.

The class of operators of the form \eqref{duality} is actually slightly larger than that given by \eqref{r-intuitive}, as it allows for the curves $\gamma_{x_2}$ and $\gamma^*_{x_1}$ to develop cusps. In general it may happen that the
kernels of $D\pi_1 (0)$ and $D\pi_2 (0)$ agree but this cannot be the case
in the setting of \eqref{r-intuitive} where $\Sigma$ is a submanifold of
$M_1 \times M_2$ and $\pi_j : \Sigma \to M_j$ are the restrictions of the
co-ordinate projections from $M_1 \times M_2$. Certain examples
of fractional integration along curves fall within the more
general class of operators given by \eqref{duality}. (See the remarks section,
Section \ref{remarks-sec}). 

\begin{definition}\label{improv-def}  Let $1 \leq p_1,p_2 \leq \infty$.  We say that the sextuple $(\Sigma, M_1, M_2, \pi_1, \pi_2, 0_\Sigma)$ is \emph{of strong-type $(p_1,p'_2)$} if we have an estimate of the form
$$ |\int_\Sigma f_1(\pi_1(x)) f_2(\pi_2(x)) a(x)\ dx| \lesssim \| f_1\|_{L^{p_1}(M_1)} \| f_2\|_{L^{p_2}(M_2)}$$
for all $f_1$, $f_2$ and all cutoff functions $a$ supported in a sufficiently small neighbourhood of $0$.  (See the Notation section for the definition of $\lesssim$).  

\end{definition}

In the intuitive setting, the above definition (modulo irrelevant technicalities when $p_2=\infty$) corresponds by duality to $R$ mapping $L^{p_1}(M_1)$ to $L^{p'_2}(M_2)$, where $p'$ is the usual dual exponent to $p$, $1/p + 1/p' = 1$.

We are interested in the problem of determining, for a fixed sextuple 
$(\Sigma, M_1, M_2, \pi_1, \pi_2, 0_\Sigma)$, the set of exponents 
$p_1, p_2$ for which the sextuple is of strong-type $(p_1, p'_2)$.  
This problem is one aspect of a much more general and difficult problem 
concerning smoothing estimates for (possibly singular or maximal) Radon 
transforms; the field is too vast to summarize here, but we refer the reader 
to the survey papers \cite{stein:kahane}, \cite{phong}, \cite{phong:icm},
\cite{gsw}, and to the recent papers \cite{cnsw}, \cite{seeger:jams}.

We will in fact work with a more convenient setting, that of restricted weak-type estimates.

\begin{definition}\label{rwt-improv-def}  Let $1 \leq p_1,p_2 \leq \infty$.  We say that the sextuple $(\Sigma, M_1, M_2, \pi_1, \pi_2, 0_\Sigma)$ is \emph{of restricted weak-type $(p_1,p'_2)$} if we have an estimate of the form
\be{e1e2}
|\int_\Sigma \chi_{E_1}(\pi_1(x)) \chi_{E_2}(\pi_2(x)) a(x)\ dx| \lesssim |E_1|^{1/p_1} |E_2|^{1/p_2}
\end{equation}
for all measurable subsets $E_1 \subset M_1$, $E_2 \subset M_2$ and all cutoff functions $a$ supported in a sufficiently small neighbourhood of $0$.  Here we use $|E|$ to denote the measure of $E$ with respect to the Riemannian metric.
\end{definition}

Note that we have the trivial bounds of $O(|E_1|)$ and $O(|E_2|)$ for the left-hand side of \eqref{e1e2} just from the hypothesis that $\pi_1$ and $\pi_2$ are local submersions.  Thus we automatically have restricted weak-type when $p'_2 \leq p_1$.  (A variant of this argument also gives strong-type in the same region).  Henceforth we restrict ourselves to the case $p_1 < p'_2$.

Clearly strong-type $(p_1,p'_2)$ implies restricted weak-type $(p_1,p'_2)$.  Conversely if $1 < p_1 < p'_2 < \infty$ and one has restricted weak-type in an open neighborhood of $(p_1,p'_2)$, then the Marcinkiewicz interpolation theorem gives strong-type $(p_1,p'_2)$.  Thus if one is willing to lose endpoints, it suffices to consider the restricted weak-type problem.

We can reformulate Definition \ref{rwt-improv-def} in a more geometric manner which is more convenient to work with.

\begin{proposition}\label{rwt-equiv}  Let $1 \leq p_1 < p'_2 \leq \infty$.  Then the sextuple is $(\Sigma, M_1, M_2, \pi_1, \pi_2, 0_\Sigma)$ is of restricted weak-type $(p_1, p'_2)$ if and only if one has the estimate
\be{om-bound}
|\Omega| \lesssim |\pi_1(\Omega)|^{1/p_1} |\pi_2(\Omega)|^{1/p_2}
\end{equation}
for all sets $\Omega \subset \Sigma$ in a sufficiently small neighborhood of 0.
\end{proposition}

\begin{proof}
If \eqref{om-bound} holds, then \eqref{e1e2} follows by setting $\Omega := \pi_1^{-1}(E_1) \cap \pi_2^{-1}(E_2) \cap B$, where $B$ is a sufficiently small neighborhood of $0$.  Conversely, if \eqref{e1e2} holds, then \eqref{om-bound} follows by setting $E_1 = \pi_1(\Omega)$ and $E_2 = \pi_2(\Omega)$.
\end{proof}

Thus, the question is to determine to what extent the size of a set $\Omega$ in $\Sigma$ is controlled by the size of its two projections $\pi_1(\Omega) \subset M_1$ and $\pi_2(\Omega) \subset M_2$.  Intuitively, the answer to this question should somehow depend on how ``independent'' the projections $\pi_1$ and $\pi_2$ are; for instance, if $M_1 = M_2$ and $\pi_1 = \pi_2$ it is easy to see that there are no estimates of the form \eqref{om-bound} other than the trivial ones when $p'_2 \leq p_1$.

To make the notion of ``independence'' more precise we introduce two vector 
fields $X_1$ and $X_2$ on $\Sigma$.  For $j=1,2$, we let $X_j$ be a smooth 
vector field defined on a neighborhood of $0_\Sigma$ such that $X_j$ is never 
zero, and $X_j$ always lies parallel to the fibers of $\pi_j$, or equivalently 
that the push-forward $(\pi_j)_* X_j$ of $X_j$ is identically zero.  
(Equivalently, $X_j(x)$ is always a non-zero element of the one-dimensional 
kernel of the derivative map $D\pi_j(x): T_x \Sigma \to T_{\pi_j(x)} M_j$).  
The vector field $X_j$ is only defined up to multiplication by smooth non-zero 
scalar functions, but we will not be bothered by this freedom and just work 
with a fixed choice of $X_1$ and $X_2$.  
From the classical Picard existence theorem for ODE we can see that such vector fields are guaranteed to exist in a sufficiently small neighborhood of 0.

With these vector fields $X_1$, $X_2$ one can then define the \emph{flow maps} $e^{tX_j}: \Sigma \to \Sigma$ in a neighborhood of 0 for sufficiently small $t$ and $j=1,2$ by the ODE
$$ \frac{d}{dt} e^{tX_j}(x) = X_j(e^{tX_j}(x)); \quad e^{0X_j}(x) = x.$$
These maps are smooth and form a group in a neighborhood of 0, i.e. $e^{sX_j} e^{tX_j} = e^{(s+t)X_j}$.  We can interpret the $X_j$ as first-order differentiation operators by the formula
$$ X_j f(x) := \frac{d}{dt} f(e^{tX_j} x)|_{t=0}.$$
The set of  first-order differential operators with smooth co-efficients is closed under Lie brackets. We denote $[X_1,X_2]$ by $X_{12}$, and observe the formula
\be{bracket}
X_{12}(x) := \frac{d^2}{dt_1 dt_2} e^{-t_1 X_1} e^{-t_2 X_2} e^{t_1 X_1} e^{t_2 X_2}(x)|_{t_1=t_2=0}.
\end{equation}
More generally, we define vector fields $X_w$ for all words $w$:

\begin{definition}\label{word-def}  We define a \emph{word} $w$ to be any non-empty finite ordered collection of 1s and 2s.  We define the \emph{degree} $\deg(w) \in \N \times \N$ to be the ordered pair $(\deg(w)_1, \deg(w)_2)$, where $\deg(w)_j$ is the number of occurrences of $j$ in $w$.  We define the vector fields $X_w$ for $w \in W$ recursively by $X_1 := X_1$, $X_2 := X_2$, and
$$ X_{wj} := [X_w, X_j] \hbox{ for all } j=1,2 \hbox{ and } w \in W,$$
thus for instance
$$ X_{12211} = [[[[X_1,X_2],X_2],X_1],X_1].$$
We give a partial ordering on degrees by writing $(a_1, a_2) \leq (b_1,b_2)$ when $a_1 \leq b_1$ and $a_2 \leq b_2$.
\end{definition}

Note that there exists a neighborhood of $0$ for which all the vector fields $X_w$ exist and are smooth.  At first glance it may appear that these vector fields $X_w$ do not cover all possible Lie bracket combinations of $X_1$ and $X_2$, but the Jacobi identity allows one to write any other Lie bracket combination as a linear combination of the $X_w$.

We say that $X_1$ and $X_2$ obey the \emph{H\"ormander condition} at 0 if 
there exist words $w_1, \ldots, w_n \in W$ such that $X_{w_1}(0), \ldots, 
X_{w_n}(0)$ span the tangent space $T_0 \Sigma$ at 0.  This condition has 
appeared several times before in work on Radon-like transforms, most notably 
in \cite{cnsw} (see also \cite{seeger:jams}) where the following
result is established. 

\begin{proposition}\label{horm}  If $X_1$ and $X_2$ do not obey the H\"ormander condition at 0, then $(\Sigma, M_1, M_2, \pi_1, \pi_2, 0)$ is not of restricted weak-type $(p_1,p'_2)$ for any $1 \leq p_1 < p'_2 \leq \infty$.  In particular, there are no non-trivial $L^p$ mapping properties near 0, either of strong type or of restricted weak-type.
\end{proposition}

For completeness we give a quick proof of this result in Section 
\ref{proof-sec}, as a consequence of the machinery developed in \cite{cnsw}.

It remains to consider the case when $X_1$ and $X_2$ do obey the H\"ormander condition, and we shall assume this for the rest of the Introduction.

\begin{definition}\label{ntuple-def}  We define $W^n$ to be the space of all $n$-tuples $(w_1, \ldots, w_n)$ of words.  If $I = (w_1, \ldots, w_n)$ is an $n$-tuple, we define the \emph{degree} $\deg(I) \in \N \times \N$ by
$$ \deg(I) = \deg(w_1) + \ldots + \deg(w_n)$$
and define the function $\lambda_I: \Sigma \to \R$ by
$$ \lambda_I(x) := \det(X_{w_1}(x), \ldots, X_{w_n}(x)).$$
\end{definition}

There exists a neighbourhood of 0 where the functions $\lambda_I$ are all well-defined smooth functions. The H\"ormander condition thus asserts that there exists an $n$-tuple $I_0$ such that $\lambda_{I_0}(0) \neq 0$.  Henceforth we fix this $n$-tuple $I_0$.

To relate these $n$-tuples $I$ to our problem \eqref{om-bound} we introduce two-parameter Carnot-Carath\'eodory balls, in the spirit of \cite{nsw}.

\begin{definition}\label{carnot-def}  Let $x \in \Sigma$ be in a sufficiently small neighbourhood of 0, and let $0 < \delta_1, \delta_2 \ll 1$ be sufficiently small numbers.  We define the \emph{two-parameter Carnot-Carath\'eodory ball}\footnote{We chose this definition for the introduction as it is the most intuitive to visualize.  However, in our rigorous argument we shall avoid these balls, and work with an equivalent family of balls defined by exponentiating various weighted commutators of the $X_i$.  See \cite{nsw} for a detailed comparison of the two types of balls.} $B(x; \delta_1, \delta_2)$ to be the closure of the set
$$ \{ e^{t_1 \delta_1 X_1} \ldots e^{t_k \delta_k X_k}(x): k \geq 0; |t_1| + \ldots + |t_k| \leq 1 \}$$
where we adopt the periodic convention that $\delta_{j+2} := \delta_j$ and $X_{j+2} := X_j$.
\end{definition}

Informally, the ball $B(x; \delta_1, \delta_2)$ represents (up to constants) the set of all points which can be reached from $x$ by flowing for an amount $\delta_1$ in the $X_1$ direction, and $\delta_2$ in the $X_2$ direction.  From \eqref{bracket} we heuristically expect that we can also flow by $\delta_1 \delta_2$ in the $X_{12}$ direction and still stay inside (a constant dilate of) this ball; more generally, we expect to flow by $\delta^{\deg(w)}$ in the $X_w$ direction, where we adopt the notation
$$ \delta^{(n_1,n_2)} := \delta_1^{n_1} \delta_2^{n_2}.$$
Because of this, we heuristically expect $B(x;\delta_1,\delta_2)$ to look something like the convex hull of the points $x \pm \delta^{\deg(w)} X_w$.  Following these heuristics, we will eventually be able to obtain a volume estimate (inspired by a similar formula in \cite{nsw}) which is roughly of the form
\be{heur}
|B(x;\delta_1, \delta_2)| \sim \sum_I \delta^{\deg(I)} |\lambda_I(x)|
\end{equation}
where $I$ ranges over a finite set depending on $I_0$ and on certain restrictions on $\delta_1$ and $\delta_2$.

The balls $B(x;\delta_1,\delta_2)$ form excellent examples of sets $\Omega$ to test \eqref{om-bound} with.  Indeed, one heuristically expects
$$ |\pi_j(B(x;\delta_1,\delta_2))| \sim \frac{|B(x;\delta_1,\delta_2)|}{\delta_j}$$
for $j=1,2$.  Inserting this into \eqref{om-bound} with $\Omega := B(x; \delta_1, \delta_2)$, we obtain the inequality\footnote{Equivalently, one can obtain this bound by applying the original averaging operator $R$ to the characteristic function of $\pi_1(B(x;\delta_1,\delta_2))$.  These sets (which tend to look like squashed neighborhoods of curve arcs) are well-known as test sets for these averaging operators; see e.g. \cite{ccc}, \cite{gs:98}, \cite{seeger:jams}.}
$$ |B(x; \delta_1, \delta_2)| \gtrsim \delta_1^{c_1} \delta_2^{c_2}$$
as a special case of \eqref{om-bound}, where the positive numbers $c_1$, $c_2$ are given by
\be{c-def}
c_1 := p_2 / (p_1 + p_2 - p_1 p_2), \ \ c_2 := p_1 / (p_1 + p_2 - p_1 p_2).
\end{equation}  
Comparing this with \eqref{heur}, we are thus led to the following Conjecture.

\begin{definition}\label{newton-def}  We define the \emph{Newton polytope} $P$ of $X_1$ and $X_2$ at 0 to be the closed convex hull of the set
$$ \{ x \in \R^+ \times \R^+: x \geq \deg(I), \lambda_I(0) \neq 0 \hbox{ for some } I \in W^n \}.$$
\end{definition}

Note that $P$ is indeed the Newton polytope of the right-hand side of \eqref{heur}, if one thinks of this as a Taylor series in $\delta_1$ and $\delta_2$.

\begin{conjecture}\label{conj}  Let $(\Sigma, M_1, M_2, \pi_1, \pi_2, 0)$ and $X_1$, $X_2$ obey the H\"ormander condition, and let $1 \leq p_1 < p'_2 \leq \infty$.  Then we have strong type $(p_1, p'_2)$ if and only if the point $(c_1, c_2)$ defined by \eqref{c-def} lies in the closed polytope $P$.
\end{conjecture}

The main result of this paper is that this conjecture is true away from endpoint cases.

\begin{theorem}\label{main}  Let the notation be as in Conjecture \ref{conj}.  Then we have strong type $(p_1, p'_2)$ when $(c_1, c_2)$ lies in the interior of $P$.  Conversely, if $(c_1, c_2)$ lies in the exterior of $P$, then one does not even have restricted weak-type $(p_1, p'_2)$.
\end{theorem}

As one can verify from elementary algebra, the statement that $(c_1, c_2)$ lies in $P$ is equivalent to $(1/p_1, 1/p'_2)$ lying in the closed convex hull of
\be{hull}
\{ (0,0), (1,1) \} \cup \{ (\frac{\deg(I)_1}{\deg(I)_1 + \deg(I)_2 - 1},
\frac{\deg(I)_1-1}{\deg(I)_1 + \deg(I)_2 - 1}): I \in W^n, \lambda_I(0) \neq 0 \}.
\end{equation}
Similarly, the statement that $(c_1, c_2)$ lies in the interior of $P$ is equivalent to $(1/p_1, 1/p'_2)$ lying in the interior of the convex hull of \eqref{hull}.

For curves in the plane ($n=3$), Seeger \cite{seeger:jams} has obtained
Theorem 1.10 previously for the operators $R$ in \eqref{r-intuitive} (here
$\Sigma$ is a submanifold of $M_1 \times M_2$ and the $\pi_j$ 's are the
induced co-ordinate projections from $M_1 \times M_2$). In this case,
$X_1$ and $X_2$ are necessarily linearly independent at $0$ and so the
3-tuples $I\in W^3$ with $\lambda_I (0) \ne 0$ defining \eqref{hull}
reduce to looking at iterated commutators $X_w , w\in W$, which do not
lie in the subspace spanned by $X_1$ and $X_2$ at $0$. See \cite{seeger:jams}
for details.
  
We now give some examples of Theorem \ref{main}.

{\bf Convolution with curves.}  We consider operators $A_n$ mapping functions on $\R^{n-1}$ to $\R^{n-1}$ defined by
$$ A_n f(x) = \int f(x-\gamma(t)) a(t)\ dt$$
where $a$ is a smooth cutoff in a small neighborhood of 0 and $\gamma$ is a smooth curve in $\R^{n-1}$ with $\gamma(0)=0$ and $\gamma'(0) \neq 0$.  This operator is of the above type with
$$ M_1 := M_2 := \R^{n-1}; \quad \Sigma := \{ (x,t): x \in \R^{n-1}, t \in \R\}; \quad 0_{M_1} := 0_{M_2} := 0_\Sigma := 0$$
and
$$ \pi_1(x,t) := x; \quad \pi_2(x,t) := x + \gamma(t).$$
We can select the vector fields $X_1$, $X_2$ as
$$ X_1 := \partial_t; \quad X_2 := \partial_t - \gamma'(t) \cdot \nabla_x.$$
A routine calculation shows that
$$ X_{12} = -\gamma''(t) \cdot \nabla_x$$
and more generally that
$$ X_{w_1 w_2 \ldots w_k} = \epsilon \gamma^{(k)}(t) \cdot \nabla_x$$
for any word $w_1 \ldots w_k$, where $\epsilon = -1$ if $w_1 w_2 = 12$, $\epsilon = +1 $ if $w_1 w_2 = 21$, and $\epsilon = 0$ otherwise.  

Let us now specialize to the maximally curved case, when
$$ \gamma(0), \gamma'(0), \ldots, \gamma^{(n-1)}(0) \hbox{ span } \R^{n-1}.$$
This for instance is the case when
\be{poly}
\gamma(t) := (t, t^2, \ldots, t^{n-1}).
\end{equation}
Then we see that the vector fields
$$ X_1, X_2, X_{12}, X_{12w_3}, \ldots, X_{12w_3 \ldots w_{n-1}}$$
span $\R^{n-1} \times \R$ at zero, for arbitrary $w_3 \ldots w_{n-1} \in \{1,2\}$.  Furthermore, any other $n$-tuples of vector fields $X_w$ which span have equal or larger degree (with respect to the partial ordering $\geq$ on degrees).  A simple computation then shows that the convex hull of \eqref{hull} is the trapezoid with vertices 
$$(0, 0), (\frac{n^2-3n+4}{n^2-n}, \frac{n-2}{n}), (\frac{2}{n}, \frac{2n-4}{n^2-n}), (1, 1).$$
Thus our theorem gives strong-type $(L^{p_1},L^{p'_2})$ in the interior of this trapezoid, and failure of restricted type $(L^{p_1},L^{p'_2})$ outside this trapezoid (when $1 \leq p_1 < p'_2 \leq \infty$).  This result was already obtained by Christ \cite{ccc} (with earlier work in lower dimensions by 
\cite{littman}, \cite{ober:conv}, \cite{ober:conv2}, \cite{ober:conv3}, 
\cite{gs:98}, 
\cite{gsw-2}),        
at least in the polynomial case \eqref{poly}; in fact these papers go further and prove various strong-type or restricted-type estimates on the boundary.

{\bf Restricted x-ray transforms.}  We consider operators of the form
$$ R_n f(x', s) := \int f(x' + t\gamma(s),t) a(x',s,t)\ dt$$
where $n \geq 3$, $x' \in \R^{n-2}, \gamma(s) = (s,\ldots, s^{n-2}), s, t 
\in \R$ and $a(x',s,t)$ is a bump 
function.  This is the X-ray transform on $\R^{n-1}$ restricted to the line 
complex of lines whose direction lies in the curve $\{ (\gamma(s), 1):
s \in \R \}$.  This is a model example of operators studied for instance in
\cite{gelfand-1}, \cite{guillemin-2}, \cite{gu:89}, \cite{gu:90}, \cite{gu:91},
\cite{ce}, \cite{gs:98},  \cite{gsw-2}, \cite{ober:x-ray},
and elsewhere.  This operator can be placed into the above setting with $M_1 
:= \{ (x', t) \in \R^{n-1} \}$, $M_2 := \{ (x', s \} \in \R^{n-1}\}$, and
$$ \Sigma := \{ (x', t, s) \in \R^n \}$$
with 0 being the usual origin and the projection maps $\pi_1$, $\pi_2$ 
defined by
$$ \pi_1(x', t, s) := (x' + t\gamma(s), t); \quad
\pi_2(x', t, s) := (x', s).$$
The associated vector fields $X_1$, $X_2$ can thus be chosen to be
$$ X_1 = \partial_s - t \gamma'(s) \cdot \nabla_{x'}; \quad X_2 = \partial_t.$$
A computation then shows that
\bas
X_{12} = - X_{21} &= \gamma'(s) \cdot \nabla_{x'}\\
X_{121} = - X_{211} &= - \gamma''(s) \cdot \nabla_{x'} \\
X_{1211} = - X_{2111} &= \gamma'''(s) \cdot \nabla_{x'}
\end{align*}
and so forth.  Furthermore, all other commutators are zero.  Thus, we essentially only have one spanning set of vector fields:
$$ X_1, X_2, X_{12}, X_{121}, \ldots, X_{121\ldots1},$$
where the last vector field has $n-3$ consecutive ones.  The total degree of this $n$-tuple can then be computed to be
$$ (\frac{n^2-3n+4}{2}, n-1).$$
The convex hull of \eqref{hull} is thus the triangle with vertices
$$ (0,0), (1,1), ( \frac{n^2-3n+4}{n(n-1)}, \frac{n-2}{n} ),$$
and so our Theorem gives $L^p \to L^q$ mapping properties in the interior of this triangle and failure of $L^p \to L^q$ outside this triangle.  This agrees with previous results in \cite{gs:94}, \cite{gs:98}, \cite{gsw-2}, \cite{ober:x-ray}, \cite{ce}, although the results there include certain endpoints which are not obtained by our methods.

{\bf Folds and cusps.} When $\Sigma \subset M_1 \times M_2$ is a submanifold
and the $\pi_j$'s are the induced co-ordinate projections, the underlying
geometry is sometimes best expressed in terms of the microlocal picture
$$C$$
$$ {\tilde \pi_1} \swarrow \ \ \ \searrow {\tilde \pi_2}$$
$$T^* M_1 \setminus 0 \ \ \ \ \ \ \ T^* M_2 \setminus 0$$
where $C = N^* \Sigma' = \{ (x_1, \xi_1; x_2, \xi_2) \in T^* (M_1 \times M_2)
\setminus 0 : (\xi_1, -\xi_2) \perp T_{(x_1,x_2)} \Sigma \}$, the twisted
conormal bundle of $\Sigma$, is the canonical relation in $T^* M_1 \times
T^* M_2$ (we are using perhaps nonstandard notation here).
When $n\ge 4$, there must be points in $C$ where $D{\tilde \pi_2}$ and
$D{\tilde \pi_1}$ drop rank. 
The simplest singularity that can occur is a Whitney
fold ($S_{1,0}$ in the Thom-Boardman description) and the second simplest
stable singularity for ${\tilde \pi_1}$ or ${\tilde \pi_2}$, 
after folds, is a Whitney cusp;
singularity class $S_{1,1,0} := S_{1_2,0}$.        

For the restricted x-ray transform $R_n$ above, ${\tilde \pi_1}$ 
has only singularities
that belong to the Morin singularity classes $S_{1_k,0}, 1\le k \le n-3$ (see 
\cite{gs:escorial} for a discussion of these singularity classes and
their connections with the operators considered here), while
both projections ${\tilde \pi_1}$ and ${\tilde \pi_2}$ 
for the convolution example $A_n$
(when $\gamma$ is maximally curved) have singularities at most of type
$S_{1_k,0}, 1\le k \le n-3$. We refer to situations as {\it one-sided}
when conditions are imposed on one projection while no assumptions are
imposed on the other projection. On the other hand, situations 
for which both ${\tilde \pi_1}$
and ${\tilde \pi_2}$ belong to a given singularity class are referred to as
{\it two-sided}. Thus the geometry underlying $A_4$ is a two-sided fold
while that of $A_5$ is a two-sided cusp. See \cite{guillemin-2} for the
microlocal analysis of the restricted x-ray transform. 
Also see \cite{gs:94}, \cite{gs:98},
\cite{gsw}, \cite{phongstein-old} for more details. 

To see examples of how this can be expressed in our setting, consider
operators of the form
\be{diagonal}
Rf(x) = \int f(\gamma(x,t)) a(x,t) \, dt
\end{equation}
where $n\ge 3, \gamma : \R^{n-1} \times \R \to \R^{n-1}$ is a smooth
map with $\gamma(x,t) \equiv x$ and $a(x,t)$ is a bump function. Thus
$x \mapsto \gamma_t (x) := \gamma(x,t)$ is a family of diffeomorphisms 
for small $t$ and a result in \cite{cnsw} associates to $\{\gamma_t \}$
a unique sequence of smooth vector fields $Z_1, Z_2, \ldots $ 
such that $\gamma_t (x) \sim \exp (\sum t^j Z_j ) (x)$ to infinite order
as $t\rightarrow 0$. This operator can be put into the above setting
with $M_1 = M_2 = \R^{n-1}, \Sigma = \{ (x,t) : x\in \R^{n-1}, t\in \R \}$,
with $0$ being the usual origin and 
$$\pi_1 (x,t) = \gamma_t (x), \ \ \ \pi_2 (x,t) = x.$$
(Equivalently, if $\frac{\partial\gamma}{\partial t} \ne 0$, 
one could take $\Sigma = \{ (\gamma_t (x),x) : x\in \R^{n-1},
t\in \R \}$ so that $\Sigma \subset M_1 \times M_2$ is a submanifold
containing the diagonal.)
The associated vector fields $X_1$ and $X_2$ can be chosen to be
$$X_1 = W(x,t) \cdot \nabla_x + \partial_t: \ \ X_2 = \partial_t$$
where $W(x,t) = \frac{d}{dh} \gamma_{t+h}^{-1} \circ \gamma_t (x) |_{h=0}$.
Observe that $X_{12} = [X_1,X_2] = W'(x,t)\cdot \nabla_x$ where $'$
denotes differentiation with respect to $t$. Similarly
$X_{122\dots 2} = (-1)^k W^{(k)} (x,t) \cdot \nabla_x$ where $k$ is the number
of 2's.  

In this setting, when $n=4$, a result of Phong and Stein 
\cite{phongstein-old} states that ${\tilde \pi_2}$ (resp. 
${\tilde \pi_1}$) has at most
Whitney fold singularities iff $Z_1, Z_2,$ and $Z_3 + \frac{1}{6} [Z_1,
Z_2]$ (resp. $Z_1, Z_2$ and $Z_3 - \frac{1}{6} [Z_1,Z_2]$) are 
linearly independent (A similar result holds for cusps, see \cite{gs:98}). 
A simple computation, using the Baker-Campbell-Hausdorff 
formula, shows that $Z_1, Z_2$ and $Z_3 + \frac{1}{6}[Z_1,Z_2]$ are linearly
independent at $x$ iff $W(x,0), W'(x.0)$ and $W''(x,0)$ are linearly independent
and this is equivalent to 
$$\lambda_I (x,0) := {\rm det}(X_1,X_2,X_{12},X_{122})(x,0) \ne 0.$$
(More generally, for any $n\ge 4$, a theorem of Greenleaf and Seeger
\cite{gs:escorial} gives the equivalence between the linear independence
of the $\{W^{(k)}\}_{0\le k\le n-2}$, the linear independence of some
combination of the $Z_j$'s, and ${\tilde \pi_2}$ having only (strong)
Morin singularities, $S^{+}_{1_k ,0}, \ k\le n-3$.)
Since $I=(1,2,12,122)$ has degree $(3,4)$, Theorem 1.10 implies that when $n=4$
and ${\tilde \pi_2}$ has at most Whitney fold singularities, 
the operator $R$ maps $L^{p_1}$
to $L^{p_2 '}$ whenever $(1/p_1,1/p_2 ')$ 
lies in the interior of
the triangle with vertices $(0,0)$, $(1,1)$, $(1/2,1/3)$. This
is a particular case of a theorem of Greenleaf and Seeger \cite{gs:94} which 
requires an additional hypothesis but contains estimates on the boundary of the triangle. 
A similar result holds when ${\tilde \pi_1}$ has at most 
Whitney fold singularities. 

The restricted x-ray transform provides
an example of the sharpness of this one-sided fold theorem. Another
interesting example is the following one considered by
Secco \cite{secco}. 
Let 
$$\gamma_t (x) = x\cdot (t,t^2,a t^3)$$
where $\cdot$ denotes the
3-dimensional Heisenberg group multiplication. Then $W(x,t) =
(-1,-2t,-(3a+1/2)t^2 + 1/2 x_2 - x_1 t)$ and so
$$ X_{12} = 2\partial_{x_2} + ((6a+1)t + x_1)\partial_{x_3},$$
$$X_{122} =
-(6a+1)\partial_{x_3}, \ \ \ \ X_{121} = (1-6a)\partial_{x_3}.$$
All other iterated commutators are zero. When $a=1/6$, Theorem
1.10 shows that the corresponding averaging operator \eqref{diagonal}
does not map
$L^{p_1}$ to $L^{p_2 '}$ whenever $(1/{p_1},1/{p_2 '})$
lies outside the triangle with vertices $(0,0)$, $(1,1)$, $(\frac{1}{2},
\frac{1}{3})$. See \cite{secco} for these results including endpoint
estimates. 

{\bf A nondegenerate example $\R^5$.} As the final example, we consider
an operator of the form \eqref{diagonal} when $n=5$ and where 
$$\gamma(x,t) = (x_1 + t, x_2 + t^2, x_3 + t^3, x_4 + t^4 + x_2 t).$$
One easily computes that $W(x,t) = (-1,-2t,-3t^2,-4t^3 - x_2 + 2t^2)$
and so at the origin we have
$$X_1 = -\partial_{x_1}+\partial_t, \ \ X_2=\partial_t, \ \ 
X_{12}=2\partial_{x_2},$$
$$X_{121}=-6\partial_{x_3}+2\partial_{x_4}, \ \ \ X_{122}=-6\partial_{x_3} +
4\partial_{x_4},$$
$$X_{1211} = X_{1212} = X_{1221} = X_{1222} = 24\partial_{x_4};$$
all other commutators being zero. Therefore the only 5-tuples of 
iterated commutators that contribute 
have degrees $(5,5), (7,4), (6,5), (5,6)$ and $(4,7)$. The convex
hull of \eqref{hull} is then the pentagon with vertices
$$(0.0), (\frac{4}{10},\frac{3}{10}), (\frac{5}{9},\frac{4}{9}), (\frac{7}
{10}, \frac{6}{10}), (1,1).$$
Our Theorem thus gives strong type $(L^{p_1},L^{p_2 '})$ in the
interior of this pentagon, and failure of restricted weak type
$(L^{p_1},L^{p_2 '})$ outside. This pentagon is a larger region
in $(\frac{1}{p_1},\frac{1}{p_2 '})$ space than the corresponding
one for the operator $A_5$ considered above. 

In the remarks section (Section \ref{remarks-sec}) we explain 
the difficulties in
extending our arguments for Theorem 1.10 to the endpoint cases.  
After some geometric preliminaries
in Sections \ref{notation-sec}, \ref{balls-sec}, we will prove the first part
of Theorem \ref{main} in Section \ref{necessary-sec}, and the second part 
in Sections \ref{reduce-sec}-\ref{conclude-sec}. 

We would like to thank Mike Christ, Allan Greenleaf, Andreas Seeger, and the referee 
for several useful discussions and comments.  The authors also thank Mike Greenblatt for pointing out an error in the first version of this manuscript.

This work was conducted at UNSW.  The first author is a Clay Prize Fellow 
and is supported by a grant and Packard Foundations. The second author is
supported in part by an ARC grant.

\section{Proof of Proposition \ref{horm}}\label{proof-sec}

We first give a proof of Proposition \ref{horm}.  The idea is to use a quantitative version of the proof of Frobenius's theorem given by Theorem 8.8 in \cite{cnsw}.  In our language, this theorem asserts (among other things) that if $X_1$ and $X_2$ do not obey the H\"ormander condition at 0, then for any $N > 0$ there exists a submanifold $S \subset \Sigma$ containing $0$ of positive codimension such that for all $0 < \delta \ll 1$ and all $x \in S$, $t \in \R$ with $\dist(x,0) \lesssim \delta$ and $|t| \lesssim \delta$ we have
\be{flow}
\dist(e^{t X_1} x, S), \dist(e^{t X_2} x, S) \lesssim \delta^N.
\end{equation}
In other words, $S$ is ``almost invariant'' under the flow of $X_1$ and $X_2$.
The proof of this claim in \cite{cnsw} uses a quantitative form of the Baker-Campbell-Hausdorff formula; it is also possible to proceed using the machinery in this paper but we will not do so here.

Assuming this statement, we test \eqref{om-bound} with the following neighbourhood of $S$:
$$ \Omega := \{ x \in \Sigma: \dist(x,0) \lesssim \delta; \dist(x,S) \lesssim \delta^N \}.$$
If $S$ has codimension $k$, then we have
$$ |\Omega| \sim \delta^{n-k} \delta^{Nk}.$$
On the other hand, from \eqref{flow} and the smoothness of $X_1$, $X_2$ we see that the set
$$ \{ e^{tX_1} x: x \in \Omega; |t| \lesssim \delta \}$$
is contained in a set which is essentially the same as $\Omega$ (but with different implicit constants in the definition) and so also has measure $\sim \delta^{n-k} \delta^{Nk}$.  From Fubini's theorem and the fact that $\pi_1$ is a submersion we thus have
$$ |\pi_1(\Omega)| \lesssim |\Omega|/\delta.$$
Similarly we have
$$ |\pi_2(\Omega)| \lesssim |\Omega|/\delta.$$
Inserting these facts into \eqref{om-bound} we obtain the condition
$$ \delta^{n-k} \delta^{Nk} \lesssim (\delta^{n-k} \delta^{Nk} \delta^{-1})^{1/p_1 + 1/p_2};$$
since $\delta$ can be arbitrarily small, this forces
$$ (n-k+Nk) \geq (n-k+Nk-1)(\frac{1}{p_1}+\frac{1}{p_2}).$$
Since $N$ can be arbitrarily large and $k > 0$, we thus obtain
$$ \frac{1}{p_1}+\frac{1}{p_2} \leq 1$$
and the claim follows.

\section{Notation}\label{notation-sec}

In our argument the sextuplet $(\Sigma, M_1, M_2, \pi_1, \pi_2, 0)$ shall be fixed.  We assume that $X_1$, $X_2$ obey the H\"ormander condition.  By continuity we may find a fixed $n$-tuple $I_0$ and a fixed open neighbourhood $V \subset \Sigma$ of $0$ such that the vector fields $X_w$ are all smooth and well-defined on $V$, and $\lambda_{I_0}$ is bounded away from zero on $V$.
Unless otherwise specified, the symbol $x$ will always be assumed to denote a point in $V$.

In our argument we will need a large parameter $N \gg 1$, which will eventually depend on $I_0$, $p_1$, $p_2$, and we will also need a small parameter $0 < \eps \ll 1$, which will eventually depend on $N$ and $I_0$, $p_1$, $p_2$.  For the first part of the paper, in which we develop the theory of Carnot-Carath\'eodory balls, we shall only use the parameter $\eps$; the $N$ parameter will be used later when we begin to prove \eqref{om-bound}.

We use $C$ to denote various large numbers which depend on $I_0$, $p_1$, $p_2$, $V$, and the sextuplet $(\Sigma, M_1, M_2, \pi_1, \pi_2, 0)$. We use $C_N$ to denote various large numbers which depend on the above parameters as well as $N$, and $C_{N,\eps}$ to denote various large numbers which depend on the above parameters as well as $N$, and $\eps$.  We will always assume implicitly that these constants have been chosen sufficiently large.

Let $A, B$ be positive quantities.
We use $A \lesssim B$ or $A = O(B)$ to denote the estimate $A \leq C_{N,\eps} B$, and $A \ll B$ to denote the estimate $A \leq C_{N,\eps}^{-1} B$.  We use $A \sim B$ to denote the estimate $A \lesssim B \lesssim A$.  Later on we will introduce a variant notation $A \lessapprox B$.

Thus for instance we have
\be{lam-base}
\lambda_{I_0}(x) \sim 1 \hbox{ for all } x \in V.
\end{equation}

We also let $K$ be a large constant (depending only on $I_0$, $p_1$, $p_2$, $N$, $\eps$, and the sextuplet) to be chosen later; our constants $C$, $C_N$, $C_{N,\eps}$ will \emph{not} be allowed to depend on $K$ unless explicitly subscripted.

\section{Two-parameter Carnot-Carath\'eodory balls}\label{balls-sec}

We now set up the machinery needed to study the Carnot-Carath\'eodory balls.  Our arguments here are very much inspired by the beautiful paper of \cite{nsw}.  However the results in \cite{nsw} are phrased for one-parameter balls and in order to exploit them one would have to reproduce many of the arguments for the two-parameter setting.  This would have added no new insights beyond those already in \cite{nsw}, so we have adopted a different approach - based more on Gronwall's inequality than the Baker-Campbell-Hausdorff formula - to give the same type of results.

In this section we fix an $x_0 \in V$, which we assume to be sufficiently close to 0.  We also assume $\delta = (\delta_1, \delta_2)$ to be an arbitrary pair of numbers which obey the smallness condition 
$$0 < \delta_1, \delta_2 \leq C_{N,\eps,K}^{-1},$$ 
and the non-degeneracy condition
$$ \delta_1^{1/\eps} \lesssim \delta_2 \lesssim \delta_1^\eps.$$
These will be the two radii for our two-parameter Carnot-Caratheodory balls.
As in the introduction we define 
$$\delta^{(n_1,n_2)} := \delta_1^{n_1} \delta_2^{n_2}$$
and
$$
(K\delta)^{(n_1,n_2)} := (K \delta_1)^{n_1} (K \delta_2)^{n_2}
$$
for any pair $(n_1,n_2)$.

We let $\I \subset W^n$ denote the set
$$
\I := \{ I \in W^n: \deg(I) \leq \frac{2}{\eps} \deg(I_0) \}.
$$
Observe that $\I$ is a finite set containing $I_0$, and that from the
smallness and non-degeneracy conditions of $\delta_1, \delta_2$,
\be{i-def}
(K\delta)^{\deg(I)} \lesssim (K\delta)^{\deg(I_0)} \hbox{ for all } I \not \in \I.
\end{equation}

Define the vector-valued function $\Lambda$ by
\be{Lam-def}
\Lambda(x) := ((K\delta)^{\deg(I)} \lambda_I(x))_{I \in \I}.
\end{equation}
We shall see later that $|\Lambda(x)|$ controls the volume of $B(x; \delta_1, \delta_2)$, if $x$ is sufficiently close to 0.

From \eqref{lam-base} we have
\be{lam-min}
|\Lambda(x)| \gtrsim (K\delta)^{\deg(I_0)};
\end{equation}
from \eqref{i-def} and the smallness of $\delta_1$, $\delta_2$ we thus have
\be{i-bound}
(K\delta)^{\deg(I)} |\lambda_I(x)| \lesssim |\Lambda(x)|
\end{equation}
for all $I \in W^n$, with the implicit constant here allowed to depend on $I$.

Fix an $x_0 \in V$ and observe from \eqref{Lam-def} we may find an $n$-tuple $I_{x_0} = (w_1, \ldots, w_n)$ in $\I$ such that
\be{imax}
|(K\delta)^{\deg(I_{x_0})} \lambda_{I_{x_0}}(x_0)| \sim |\Lambda(x_0)|.
\end{equation}
Fix this $I_{x_0}$.  Define the map $\Phi = \Phi_{x_0}$ from $\R^n$ to $V$ by
\be{Phi-def}
\Phi(t_1, \ldots, t_n) := \exp(\sum_{j=1}^n t_j K^{-1} (K\delta)^{\deg(w_j)} X_{w_j}) x_0
\end{equation}
and consider the ``balls'' $\Phi(B)$, where $B$ is a ball in $\R^n$ centered at the origin with fixed radius $C$, which we will choose later.  It turns out that these sets are roughly comparable to the balls $B(x; \delta_1, \delta_2)$ defined in the Introduction (for a detailed discussion, see \cite{nsw}).  In fact in the proof of Theorem \ref{main} we shall work exclusively with the $\Phi(B)$ and avoid using the balls $B(x; \delta_1, \delta_2)$.  

Since $\Phi$ is smooth and
\be{deriv}
\partial_{t_j} \Phi(0) = K^{-1} (K\delta)^{\deg(w_j)} X_{w_j}
\end{equation}
we see from \eqref{imax} that 
$$ \det(D\Phi)(0) \sim K^{-n} |\Lambda(x_0)|$$
and in particular that $\Phi$ is locally invertible near the origin.  (We will make this statement more quantitative later on).  

For each $w \in W$, let $Y_w$ be the pull-back of the vector field $K^{-1} (K\delta)^{\deg(w)} X_w$ by the map $\Phi$; these are defined near 0 by the local invertibility of $\Phi$, though we shall shortly show that they are in fact defined on all of $B$ (there is a slight abuse of notation here since the $Y_w$'s
are not exactly iterated commutators of $Y_1$ and $Y_2$ but in fact a multiple
of the corresponding commutator - e.g., $[Y_1, Y_2] =
K Y_{12}$ by the way we have just defined $Y_{12}$).  We also consider the standard Euclidean vector fields $\partial_1, \ldots, \partial_n$ on $B$, as well as the co-ordinate functions $t_1, \ldots, t_n$.  From \eqref{deriv} we have
\be{deriv-pullback}
Y_{w_j}(0) = \partial_j(0).
\end{equation}
Thus the system of vector fields $Y_{w_1}, \ldots, Y_{w_n}$ is a perturbation of the Euclidean system $\partial_1, \ldots, \partial_n$ near the origin.  

We now come to the main result of this section, which gives regularity of the vector fields $Y_w$ and controls the geometry of the ball $\Phi(B)$.

\begin{proposition}\label{ball-geom}  The vector fields $Y_w$ are well-defined and smooth on all of $B$, and obey the estimates
\be{ywm-smooth}
\| Y_w \|_{C^M(B)} \lesssim C_{w,M}
\end{equation}
for all $w \in W$ and $M \geq 0$, if $K$ is sufficiently large depending on $M$. 

For the vector fields $Y_{w_i}$ we have the more precise bounds
\be{y-bound}
Y_{w_i}(t) = \partial_i + O(|t|/K)
\end{equation}
for $t \in B$ and $i = 1, \ldots, n$.  In particular we have
\be{det-large}
\det(Y_{w_1}, \ldots, Y_{w_n})(t) \sim 1 \hbox { for all } t \in B.
\end{equation}
Finally, we have the volume bounds
\be{volume}
|\Phi(E)| \sim K^{-n} |\Lambda(x_0)| |E|
\end{equation}
for any subset $E$ of $B$.
\end{proposition}

\begin{proof}

{\bf Step 1.  Estimate $Y_{w_j}$.}

The first step is to bound the $Y_{w_i}$, and in particular to prove \eqref{y-bound} and \eqref{det-large}.

By \eqref{Phi-def} and the definition of exponentiation we have the vector field identity
\be{key-ident}
\sum_{j=1}^n t_j \partial_j = \sum_{j=1}^n t_j Y_{w_j}
\end{equation}
on $B$.  

Let $1 \leq i \leq n$.  Taking the Lie bracket of \eqref{key-ident} with $Y_{w_i}$, we obtain
$$
\sum_{j=1}^n Y_{w_i}(t_j) \partial_j + t_j [Y_{w_i}, \partial_j]
= \sum_{j=1}^n Y_{w_i}(t_j) Y_{w_j} + t_j [Y_{w_i}, Y_{w_j}]$$
which we rewrite as
$$ \sum_{j=1}^n t_j [\partial_j, Y_{w_i} - \partial_i] + (Y_{w_i} - \partial_i)
= -\sum_{j=1}^n (Y_{w_i} - \partial_i)(t_j) (Y_{w_j}-\partial_j)
- \sum_{j=1}^n t_j [Y_{w_i}, Y_{w_j}].$$
If we write
\be{y-decomp}
Y_{w_i} =: \partial_i + \sum_{k=1}^n a_i^k(t) \partial_k
\end{equation}
then the previous becomes
\be{diff}
\sum_{j=1}^n \sum_{k=1}^n t_j (\partial_j a_i^k) \partial_k + \sum_{k=1}^n a_i^k \partial_k
= - \sum_{j=1}^n \sum_{k=1}^n a_i^j a_j^k \partial_k - \sum_{j=1}^n t_j [Y_{w_i}, Y_{w_j}].
\end{equation}
Define the vector valued quantity
$$
A(t) := (a_i^j(t))_{1 \leq i,j \leq n};
$$
from \eqref{y-decomp}, \eqref{deriv-pullback} we observe that $A(0)=0$.  We now try to use \eqref{diff} to obtain good bounds on $A$.

From \eqref{decomp} below we see that $[Y_{w_i}, Y_{w_j}]$ can be written as a linear combination of the $Y_{w_k}$ with co-efficients $O(K^{-1})$, so that the $\partial_k$ co-efficient of this vector field is $O((1 + |A|)/K)$.  From \eqref{diff} we thus have
$$
|\sum_{j=1}^n t_j (\partial_j a_i^k) + a_i^k| \lesssim |A|^2 + (r|A| + r)/K
$$
for all $k$, where $r := (\sum_{i=1}^N t_i^2)^{1/2}$ is the radial co-ordinate.  We can rewrite this as
$$ |\partial_r(r a_i^k)| \lesssim |A|^2 + (r|A| + r)/K$$
where $\partial_r = \sum_{j=1}^n \frac{t_j}{r} \partial_j$ is the radial vector field.  From the definition of $A$ we thus have
$$ | \partial_r(r A)| \lesssim |A|^2 + (r|A| + r)/K.$$
On the other hand, we have $A = 0$ when $r=0$, so we have the integral inequality
$$ r |A(r\omega)| \lesssim  \int_0^r |A(s\omega)|^2 + \frac{1}{K} (s |A(s\omega)| + s ) \ ds$$
for all $\omega$ on the unit sphere.  In particular we have
$$ \frac{K |A(r\omega)|}{r} \lesssim 1 + \frac{1}{K} (\sup_{0 \leq s \leq r} \frac{K |A(s\omega)|}{s} + (\sup_{0 \leq s \leq r} \frac{K |A(s\omega)|}{s})^2)$$
for $r = O(1)$.  If $K$ is sufficiently large, we thus see from standard continuity arguments (recall that $A$ is smooth and vanishes at the origin) that
$A$ is in fact defined for all $r$ and that
$\frac{K|A(r\omega)|}{r} \lesssim 1$
for all $r\omega \in B$.  In other words we have
\be{wash}
A(t) = O(r/K)
\end{equation}
for all $t \in B$, and \eqref{y-bound} and \eqref{det-large} follows.

{\bf Step 2.  Control $\Lambda$ and $\lambda_{I_{x_0}}$.}

The next step is to understand the behaviour of the vector-valued function $\Lambda$.  We begin by proving some regularity properties of $\Lambda$ in the directions $X_w$.

Let $I \in \I$, $w \in W$, $x \in V$ be arbitrary.
We observe the identity
\be{iter}
\begin{split}
X_w \det(X_{w_1}, \ldots, X_{w_n}) = &{\rm div}(X_w) \det(X_{w_1}, \ldots, X_{w_n}) \\
&+ \sum_{j=1}^n \det(X_{w_1}, \ldots, [X_w, X_{w_j}], \ldots, X_{w_n});
\end{split}
\end{equation}
this formula is just the Lie derivative of 
$X_{w_1} \wedge \ldots \wedge X_{w_n}$
with respect to the vector field $X_w$ and  can be deduced 
(for instance) by writing everything in co-ordinates and using the product rule together with the formula
$$ ([X,Y](x))_i = \sum_{j=1}^n X_j(x) \partial_j Y_i(x) - Y_j(x) \partial_j X_i(x).$$

From \eqref{iter} and the observation (using the Jacobi identity) that $[X_w, X_{w_j}]$ is a linear combination of those $X_{w'}$ with $\deg(w') = \deg(w) + \deg(w_j)$, we see from \eqref{i-bound} that
\be{ld-1}
|(K\delta)^{\deg(w) + \deg(I)} X_w \lambda_I(x)| \lesssim C_w |\Lambda(x)|
\end{equation}
From this and \eqref{Lam-def} we obtain
\be{ld-2}
|(K\delta)^{\deg(w)} X_w \Lambda(x)| \lesssim C_w |\Lambda(x)|
\end{equation}
By iterating \eqref{iter} we thus obtain
\be{ld-3}
|(K\delta)^{\deg(w_1) + \ldots + \deg(w_k) + \deg(I)} 
X_{w_1} \ldots X_{w_k} \lambda_I(x)| \lesssim C_{w,k} |\Lambda(x)|
\end{equation}
for any $k \geq 1$; note that any term of the form ${\rm div}(X_w)$ is smooth and so will behave very nicely with respect to derivatives. 

From \eqref{ld-2} and \eqref{Phi-def} we obtain 
$$ \partial_r |\Lambda( \Phi(t) )| \lesssim K^{-1} |\Lambda( \Phi(t) )|$$
for all $t \in B$.  Similarly, from \eqref{ld-1} and \eqref{Phi-def}
we have
$$
|\partial_r (K\delta)^{\deg(I_{x_0})} \lambda_{I_{x_0}}( \Phi(t) )| 
\lesssim K^{-1} |\Lambda( \Phi(t) )|.$$
From these estimates, Gronwall's inequality, and \eqref{imax} we thus obtain
\be{lam-big}
|(K\delta)^{\deg(I_{x_0})} \lambda_{I_{x_0}}(x)| \sim |\Lambda(x)| \sim |\Lambda(x_0)|
\end{equation}
for all $x \in \Phi(B)$, if $K$ is sufficiently large.  In particular we have
\be{lam-const}
\lambda_{I_{x_0}}(x) \sim \lambda_{I_{x_0}}(x_0)
\end{equation}
for all $x \in \Phi(B)$.

{\bf Step 3.  Control $Y_w$.}

Having controlled $Y_{w_i}$, and having controlled the determinants $\lambda_I$, we will now be able to control all the other $Y_w$, and in particular $Y_1$ and $Y_2$.  As a consequence we will be able to prove \eqref{ywm-smooth}.
 
From Cramer's rule we see that for all $w \in W$ we may write
\be{decomp}
(K\delta)^{\deg(w)} X_w(x) = \sum_{j=1}^n c_{w,j}(x) (K\delta)^{\deg(w_j)} X_{w_j}(x)
\end{equation}
for all $x \in \Phi(B)$, where the co-efficients $c_{w,j}$ are linear combinations of ratios 
\be{ratios}
\frac{ (K \delta)^{\deg(I)} \lambda_I(x) }{ (K \delta)^{\deg(I_{x_0})} \lambda_{I_{x_0}}(x) }
\end{equation}
for various $I \in W^n$.  In particular, from \eqref{i-bound} we see that $c_{w,j} = O(1)$ (with the implicit constant depending on $w$).

Pulling back \eqref{decomp} under $\Phi$ we obtain
\be{yw-decomp}
Y_w = \sum_{j=1}^n \tilde c_{w,j} Y_{w_j}
\end{equation}
where $\tilde c_{w,j} := c_{w,j} \circ \Phi$. 

Fix $M \geq 0$, $1 \leq j, k \leq n$, $w \in W$ and $t \in B$.  To prove \eqref{ywm-smooth}, it suffices by \eqref{yw-decomp} to show the bounds
\be{cb}
|\nabla^M_t \tilde c_{w,j}(t)|\lesssim C_{M,w}
\end{equation}
and
\be{ab}
|\nabla^M_t a_j^k(t)| \lesssim C_{M}.
\end{equation}

The claims \eqref{cb}, \eqref{ab} have already been proven for $M=0$.  Now suppose inductively that $M > 0$ and that \eqref{cb}, \eqref{ab} have already 
been proven for all smaller values of $M$.  Recall that $c_{w,j}$ is a linear combination of ratios \eqref{ratios}.  By repeated applications of the quotient rule, \eqref{i-bound}, \eqref{lam-big}, and \eqref{ld-3} we thus have that
$$ (K \delta)^{\deg(w_{i_1}) + \ldots + \deg(w_{i_M})} |X_{w_{i_1}} \ldots X_{w_{i_M}} c_{w,j}(x)| \lesssim C_{M,w}$$
for all $w \in W$, $x \in \Phi(B)$ and $1 \leq i_1, \ldots, i_M \leq n$.  Pulling this back by $\Phi$ and throwing away some powers of $K$, we obtain
$$ |Y_{w_{i_1}} \ldots Y_{w_{i_M}} \tilde c_{w,j}(t)| \lesssim C_{M,w}$$
for all $w \in W$, $t \in B$ and $1 \leq i_1, \ldots, i_M \leq n$.  Expanding out the derivative operators $Y_{w_{i_j}}$ using \eqref{y-bound} and the induction hypothesis \eqref{cb}, \eqref{ab} for smaller values of $M$, we obtain
$$ |\partial_{i_1} \partial{i_2} \ldots \partial_{i_M} \tilde c_{w,j}(t)|
\lesssim C_{M,w} (1 + \frac{1}{K} |\nabla^M \tilde c_{w,j}|)$$
and \eqref{cb} follows if $K$ is sufficiently large depending on $M$.

Now we show \eqref{ab}.  Let $D$ be a constant co-efficient operator of order $M$.  We take the Lie bracket of $D$ with \eqref{diff}.  From the Euler identity
$$ \sum_{j=1}^n D t_j (\partial_j a_i^k) = M D a_i^k + \sum_{j=1}^n t_j (\partial_j D a_i^k)$$
we obtain
$$
\sum_{j=1}^n \sum_{k=1}^n t_j (\partial_j D a_i^k) \partial_k + (M+1) \sum_{k=1}^n (D a_i^k) \partial_k
= - \sum_{j=1}^n \sum_{k=1}^n D(a_i^j a_j^k) \partial_k - \sum_{j=1}^n [D, t_j [Y_{w_i}, Y_{w_j}]].
$$
We now estimate the terms on the right-hand side.
From the inductive hypothesis \eqref{ab}, \eqref{wash} and the Leibnitz rule we have
$$ |D(a_i^j a_j^k)| \leq C_M (1 + r |\nabla^M A|).$$
We can write $[Y_{w_i}, Y_{w_j}]$ using the Jacobi identity as a linear combination of $Y_w$.  From \eqref{yw-decomp}, the Leibnitz rule, the inductive hypotheses \eqref{cb}, \eqref{ab} for smaller values of $M$, and the claim \eqref{cb} for $M$ that was just proven, we thus have
$$ |[D, t_j [Y_{w_i}, Y_{w_j}]]| \leq C_M (1 + r |\nabla^M A|).$$
Combining these facts together, we obtain
$$ |\sum_{k=1}^n r (\partial_r D a_i^k) \partial_k + (M+1) \sum_{k=1}^n (D a_i^k) \partial_k| \leq C_M(1 + r |\nabla^M A|).$$
Taking $\partial_k$ co-efficients and then letting $i$ vary, we thus have
$$ |r \partial_r D A + (M+1) DA| \leq C_M (1 + r |\nabla^M A|).$$
Multiplying by $r^M$, and letting $D$ vary, this becomes
$$ | \partial_r( r^{M+1} \nabla^M A )| \leq C_M (r^M + r^{M+1} |\nabla^M A|).$$
By Gronwall's inequality (and noting from the a priori smoothness of all quantities that $r^{M+1} \nabla^M A$ vanishes at the origin), we thus have
$$ |r^{M+1} \nabla^M A| \leq C_M r^{M+1},$$
and \eqref{ab} follows.

{\bf Step 4.  Proof of \eqref{volume}.}

From the chain rule and the definition of the $Y_w$ we have the identity
$$ \det(D\Phi) = \frac{\det(K^{-1}(K\delta)^{\deg(w_1)} X_{w_1}, \ldots, K^{-1}(K\delta)^{\deg(w_n)} X_n)}{\det(Y_{w_1}, \ldots, Y_{w_n})}.$$
From \eqref{det-large}, \eqref{lam-const}, \eqref{imax} we thus have
$$
\det(D\Phi)(t) \sim K^{-n} |\Lambda(x_0)| \hbox{ for all } t \in B.
$$
To prove \eqref{volume} it thus suffices to show that $\Phi$ is essentially one-to-one on $B$.  More precisely, we will be able to cover the large ball $B$ by $O(1)$ small balls of radius $\sim 1$ such that $\Phi$ is one-to-one on each ball.

To prove this, suppose for contradiction that there existed $t, t'$ in $B$ with $t \neq t'$, $|t-t'| \ll 1$ and $\Phi(t) = \Phi(t')$.  From \eqref{y-bound}, \eqref{ab}, \eqref{det-large} and the inverse function theorem we see that there exists a non-zero linear combination $Y = \sum_{i=1}^N \alpha_i Y_{w_i}$ with $\alpha_i = O(1)$ such that $t' = \exp(Y) t$.  This implies that $\Phi(t') = \exp(X)\Phi(t)$, where $X := \sum_{i=1}^N \alpha_i K^{-1} (K\delta)^{\deg(w_i)} X_{w_i}$.  However from \eqref{lam-big} we see that $X$ is a smooth non-vanishing vector field, with a $C^1$ norm of $O(\delta_1 + \delta_2)$ (for instance).  It is easy to see that for such vector fields, $\exp(X)$ cannot have any fixed points; indeed for such fields we differentiate the identity
$$ \frac{d}{dt} \exp(tX)(x) = X( \exp(tX)(x) )$$
in time to obtain
$$ \frac{d^2}{dt^2} \exp(tX)(x) = O( (\delta_1 + \delta_2) |\frac{d}{dt} \exp(tX)(x)| )$$
and from Gronwall's inequality we thus have 
$$ \frac{d}{dt} \exp(tX)(x) = X(x) + O((\delta_1 + \delta_2)|t| |X(x)|)$$
for $|t| \leq 1$, and the claim follows if $\delta_1 + \delta_2$ is sufficiently small.  Thus we obtain the desired contradiction. 

\end{proof}

We remark that one can use Proposition \ref{ball-geom} to show that 
$$ |B(x_0; \delta_1, \delta_2)| \sim |\Lambda(x_0)|$$
by using \eqref{volume}, \eqref{ywm-smooth} and (iterated versions of) \eqref{bracket}; we omit the details.  In particular we can prove a rigorous version of \eqref{heur} (with $I$ restricted to $\I$).

\section{Necessity of the Newton polytope}\label{necessary-sec}

We now have enough machinery to prove the easy direction of Theorem \ref{main}.  More precisely, we assume $1 \leq p_1 < p'_2 \leq \infty$ are such that the associated point $(c_1, c_2)$ defined by \eqref{c-def} is outside of the Newton polytope $P$.

The (quarter-infinite) polytope $P$ only has a finite number of vertices, all of which are of the form $\deg(I)$ for some collection of $n$-tuples $I \in W^n$.  By choosing $\eps$ sufficiently small, one can assume that all of these $n$-tuples are in $\I$.

Since $(c_1, c_2)$ is outside of $P$, we may find a half-plane of the form
$$ \{ (x_1, x_2): a_1 x_1 + a_2 x_2 \geq 1 \}$$
which contains $P$ but does not contain $c_1, c_2$.  Since $P$ is quarter-infinite and contained in the quadrant $\{ (x_1,x_2): x_1, x_2 \geq 1 \}$ we may take $0 < a_1, a_2 < 1$; by taking $\eps$ sufficiently small we may assume that $\eps \leq a_1, a_2$ and that $a_1 c_1 + a_2 c_2 \leq 1-\eps$.  

We let $0 < \delta_0 \ll 1$ be a small parameter, and set 
$$\delta := (\delta_1, \delta_2) := (\delta_0^{a_1}, \delta_0^{a_2}).$$
Observe from construction that
$$ \delta^{(c_1,c_2)} \geq \delta_0^{-\eps} \delta^{(x_1,x_2)}$$
for all $(x_1,x_2) \in P$.  In particular we have
\be{delta-large}
\delta^{(c_1,c_2)} \gtrsim C_K^{-1} \delta_0^{-\eps} |\Lambda(0)|.
\end{equation}

We now apply the machinery of the previous section, with $x_0 := 0$ and $\delta$ set as above.  We set
$$ \Omega := \Phi(B_{1/K}),$$
where $B_{1/K}$ is the ball of radius $1/K$ centered at 0.  From \eqref{volume} we have
\be{om-size}
|\Omega| \sim C_K^{-1} |\Lambda(0)|.
\end{equation}
On the other hand, from \eqref{ywm-smooth} we see that
$$ e^{tY_j} B_{1/K} \subset B$$
for all $j = 1,2$ and $|t| \leq 1/K$.  Thus
$$ e^{t \delta_j X_j} \Omega \subset \Phi(B)$$
for all $j=1,2$ and $|t| \leq 1/K$.  In particular, from \eqref{volume} we have
$$ | \bigcup_{t: |t| \leq 1/K} e^{t \delta_j X_j} \Omega| \lesssim C_K |\Lambda(0)|.$$
Since $X_j$ lies in the kernel of $d\pi_j$, is bounded away from zero, and $\pi_j$ is a submersion, we thus have
$$ \frac{1}{K} |\delta_j \pi_j(\Omega)| \lesssim C_K |\Lambda(0)|.$$
Inserting this bound and \eqref{om-size} into \eqref{om-bound} we obtain
$$ |\Lambda(0)| \lesssim C_K (|\Lambda(0)|/\delta_1)^{1/p_1} (|\Lambda(0)|/\delta_2)^{1/p_2}$$
which after some algebra and \eqref{c-def} becomes
$$ |\Lambda(0)| \gtrsim C_K \delta^{(c_1, c_2)}.$$
But this contradicts \eqref{delta-large} if we choose $\delta_0$ sufficiently small.  This concludes the proof of necessity for $(c_1,c_2)$ to lie in the Newton polytope $P$.

\section{Sufficiency of the interior of the polytope: notation and preliminary reductions}\label{reduce-sec}

We now begin the proof of the more difficult direction of Theorem \ref{main}.  For this direction we assume $1 \leq p_1 < p'_2 \leq \infty$ is such that the exponent pair $(c_1,c_2)$ defined by \eqref{c-def} lies in the interior of the polytope $P$, and our task is to prove \eqref{om-bound} (the strong-type estimates then following by Marcinkiewicz interpolation).

We may assume that the set $V$ is contained within a $1/K$-neighborhood of the origin.  In particular, if $a$ is a Lipschitz function on $V$ such that $|a(0)| \sim 1$, then $|a(x)| \sim 1$ for all $x \in V$ (if $K$ was chosen sufficiently large).

Fix $p_1$, $p_2$, $c_1$, $c_2$, and let $\Omega \subset V$ be a fixed set of positive measure.  We define the quantities $\alpha_j = \alpha_j(\Omega)$ for $j = 1,2$ by
\be{alpha-def}
\alpha_j := \frac{|\Omega|}{|\pi_j(\Omega)|};
\end{equation}
these are the analogues of the $\delta_j$ in the previous section.  Note that
$0 < \alpha_j \lesssim 1/K$ (since $V$ has diameter $O(1/K)$).  By the same algebra used in the previous section, we can rewrite \eqref{om-bound} as
\be{targ-proto}
|\Omega| \gtrsim \alpha_1^{c_1} \alpha_2^{c_2}.
\end{equation}
We adopt the periodic notation that $\alpha_{j+2} = \alpha_j$ for all $j$.  By symmetry we may assume that 
$$ \alpha_1 \geq \alpha_2.$$

Intuitively, the proof of \eqref{targ-proto} runs as follows\footnote{Our
argument follows the iteration ideas of Christ \cite{ccc} (which in turn were inspired by the work of Bourgain, Wolff, and Schlag on Kakeya-type problems) adapted to the general vector field setting.  Of course, the idea of iterating an operator supported mostly on one-dimensional sets to increase the dimension of the support is hardly new (see e.g. \cite{christ:weak}, \cite{christ:rough}, \cite{christ:curves}, \cite{greenblatt}); our main innovation is the utilization of the concept of width introduced in the previous section.}.  From \eqref{alpha-def}, we expect that for ``most'' points $x_0$ in $\Omega$, one can flow using $e^{tX_1}$ for some random $|t| \lesssim 1$, and return to $\Omega$ with probability $\sim \alpha_1$.  Similarly for $e^{tX_2}$ and $\alpha_2$.  In particular, if one considers a generic expression of the form
\be{img}
e^{t_n X_n} \ldots e^{t_1 X_1} x_0
\end{equation}
then such a point should have a probability at least $\gtrsim \alpha_n \ldots \alpha_1$ of lying in $\Omega$.  If this expression as a function of $t = (t_1, \ldots, t_n)$ had good injectivity properties (in that its Jacobian was bounded away from zero), this would then give a significant improvement to \eqref{targ-proto}, namely
$$ |\Omega| \gtrsim \alpha_n \ldots \alpha_1 =
\alpha_1^{\lfloor (n+1)/2 \rfloor}
\alpha_2^{\lfloor n/2 \rfloor}.$$
This however is too good to be true, as the results of Section \ref{necessary-sec} already indicate.  The problem is that the Jacobian of the map
$$ (t_1, \ldots, t_n) \to e^{t_n X_n} \ldots e^{t_1 X_1} x_0$$
can degenerate to zero on a complicated set (although the H\"ormander condition does ensure that this Jacobian does not vanish entirely).

In the work of Christ\cite{ccc}, the operator of convolution with the curve $(t,t^2, \ldots, t^{n-1})$ was considered.  This curve has special algebraic properties (for instance, the Jacobian turns out to be a Vandermonde determinant), and one could obtain a large set of times $(t_1, \ldots, t_n)$ where the Jacobian was reasonably large\footnote{Strictly speaking, the argument in \cite{ccc} requires an iteration by $2n-2$ steps instead of $n$ in order to locate a sub-$n$-tuple of times which can avoid the zero set of the Jacobian.  This additional complication is avoided in our approach by restricting our times to a central set of small width; if the times are allowed to range over long distances (much greater than their natural width) then one can give counter-examples to the above iteration scheme working if one only uses $n$ iterations (Christ, personal communication).} by combinatorics.  However in our more general context the Jacobian has no special algebraic structure and the zero set is about as badly behaved as a generic algebraic variety.  While algebraic varieties do have some structure which can be exploited (for instance, their intersections with other algebraic sets have a bounded number of connected components), it is not clear how a general set of times $t_1, \ldots, t_n$ can avoid this zero set.

To get around the possible algebraic complexity, we go back and do a preliminary pruning of the set of times $t_1, \ldots, t_n$ which we use in the above argument.  In particular, by using a simple stopping time argument we can restrict each time $t_i$ to a central set of fixed width (the central nature comes because the set of times is essentially a difference set), although it will cost us an epsilon in the exponents to do so.  Once one has restricted the time sets in this manner then \eqref{img} is restricted to one of the two-parameter Carnot-Carath\'eodory balls 
constructed earlier.
By using a rescaling adapted to the vector fields (specifically, we use the map $\Phi$ constructed in Section \ref{balls-sec}) we can assume that these widths are close to 1, at which point one can iterate Lemma \ref{lower} below to give satisfactory lower bounds on the Jacobian.  (This trick was inspired by the ``two ends reduction'' of Wolff \cite{wolff:kakeya}).  The condition that $(c_1, c_2)$ lie in the interior of $P$ then comes naturally from the Jacobian of the rescaling transformation (i.e. from \eqref{volume}), by computations similar to those in Section \ref{necessary-sec}.

To make the above argument rigorous we shall need some additional notation.
Let $A \lessapprox B$ denote the estimate 
$$A \leq C_{K,N,\eps} \alpha_2^{-C_N \eps} B,$$
and write $A \approx B$ if $A \lessapprox B$ and $B \lessapprox A$.
If $E$ and $F$ are sets, we say that $E$ is a \emph{refinement} of $F$, or $E \prec F$, if $E \subset F$ and $|E| \approx |F|$.

To prove \eqref{targ-proto} it will suffice to show that
\be{targ} |\Omega| \gtrapprox \alpha_1^{c_1} \alpha_2^{c_2 + C/N}.
\end{equation}
Indeed, if \eqref{targ} held, then by making $N$ large and \emph{then} making $\eps$ small, we obtain \eqref{targ-proto} with $c_1, c_2$ replaced by small perturbations of $c_1$, $c_2$.  But since $(c_1,c_2)$ is an arbitrary \emph{interior} point of $P$, and $\alpha_1, \alpha_2 \lesssim 1$, we thus obtain \eqref{targ-proto} for all $(c_1, c_2)$ in the interior of $P$ as desired.

It remains to show \eqref{targ}.  This will occupy the remaining sections of the paper, after a digression on widths.

\section{An interlude on widths}\label{width-sec}

We now introduce the main innovation of this paper, that of a central one-dimensional set with a fixed width.

\begin{definition}  Let $S$ be a subset of $[-C,C]$ with non-zero measure.  We say that $S$ is \emph{central with width} $w$ for some $w > 0$ if 
$$ |x| \lesssim w \hbox{ for all } x \in S$$
and
\be{neg}
 |I \cap S| \lesssim (|I|/w)^\eps |S|
\end{equation}
for all intervals $I$.
\end{definition}

Thus central sets with width $w$ are spread out somewhat evenly within an interval $[-Cw,Cw]$ which is centered at the origin.  It is essential that this interval is \emph{central} (i.e. centered at the origin), and that there is absolutely no portion of $S$ outside this interval, otherwise the iteration argument we use will not give a set contained in a small Carnot-Carath\'eodory ball.

Note that central sets have diameter comparable to their width.  In particular:

\begin{corollary}\label{width-mono}
If one central set of width $w$ is a subset of another central set of width $w'$, then $w \lesssim w'$. 
\end{corollary}

A central set $S$ with width $w$ is supported on an interval of length $\sim w$ and is not concentrated on any smaller interval.  This non-concentration gives such sets good properties when it comes to obtaining lower bounds of integrals of these sets.  The idea of using such a non-concentration condition was inspired by the work of Wolff \cite{wolff:kakeya} on the Kakeya problem, in which he utilized a very similar ``two-ends'' condition to achieve a similar effect.

In a later proposition (Lemma \ref{ref-lemma}) we shall construct some central sets with some width $w$.  In the remainder of this section, we show how the width property is useful.  

\begin{lemma}\label{lower} Let $P(t)$ be a polynomial of one real variable
of degree $d = 0(1)$ and $\|P\|_{\infty} \lessapprox 1$ on an interval
$[-C,C]$.  Let $S$ be a central set in $[-C,C]$ of
some width $0 < w \ll 1$.  Suppose that we have the lower bound
\be{lower-hyp}
|(\frac{d}{dt})^j P(0)| \gtrapprox w^m
\end{equation}
for some $j \geq 0$ and $m \geq 0$.  Then, if $w$ is sufficiently small, we can find a subset $S' \subset S$ of measure
$$ |S'| \approx |S|$$
such that
\be{lower-conc}
|P(t)| \gtrapprox w^{C(m,d)} \hbox{ for all } t \in S'.
\end{equation}
\end{lemma}

Observe how the hypothesis \eqref{lower-hyp} only requires a lower bound at a single point, while the conclusion \eqref{lower-conc} yields a lower bound on a large fraction of the set $S$.  Such a bound would not be possible without the non-concentration assumption \eqref{neg} (unless one was willing to lose a power of $|S|$ in \eqref{lower-conc}, which would almost certainly lead to much worse exponents\footnote{Admittedly, if one was very careful then one might still be able to obtain near-optimal results, see e.g. \cite{ccc}.  However the general theory of lower bounds is far from optimal at this stage; see \cite{ccw} for some recent progress.} $p_1, p_2$ for the positive results in Theorem \ref{main}).
 
Because of the unspecified power of $w$ in \eqref{lower-conc} we will only be able to usefully apply this lemma when the set $S$ has width close to 1.  However
we will be able to use the  map $\Phi$ developed in the previous section to rescale sets of small width to sets of width $\approx 1$ to resolve this problem.
 
\begin{proof}
By considering the homogeneous functional
$$ \| P \| := \max_{k\le d} \inf_{t\in [-C,C]} |P^{(k)}(t)|$$
on the finite dimensional space of polynomials of degree
$d$, we
see that \eqref{lower-hyp} implies
\be{lower-global} 
|P^{(k)}(t)| \gtrapprox w^m \ {\rm on} \ [-C,C]
\end{equation}
for some $k\le d$. We may assume $k\ge 1$.
Consider the set of real numbers $t$ such that $P(t) = O(w^{m + 2d})$.  
This is contained in the union of at most $O(d)$ intervals,
each interval has length at most $0(w^2)$ by \eqref{lower-global}. 
By \eqref{neg}, each such interval thus contains at most $O(w^\eps)$ of 
the set
$S$.  If we then set $S'$ to be the portion of $S$ which is outside 
these intervals, we have $|S'| \approx |S|$ (if $w$ is sufficiently small) 
and $|P(t)| \gtrapprox
w^{m + 2d}$ on $S'$ as desired.
\end{proof}

\section{Main argument}

We now give the rigorous proof of \eqref{targ}. 
Let $\I_0 \subset \I$ denote those elements of $\I$ for which $\lambda_I(0) \neq 0$.  Since $\I_0$ is finite and non-empty, we may assume by continuity that $\lambda_I(x) \sim 1$ for all $x \in V$.  In particular we have
\be{lambda-bound}
|\Lambda(x)| \gtrsim \sum_{I \in \I_0} \delta^{\deg(I)}
\end{equation}
for all $x \in V$ and $0 < \delta_1, \delta_2 \ll 1$.

Choose $\eps > 0$ small enough so that $\I_0$ contains all the vertices of $P$. 
Since $(c_1,c_2)$ lies in the polytope $P$ and $0 < \alpha_1, \alpha_2 \lesssim 1$, we see from construction of $P$ that
$$\alpha_1^{c_1} \alpha_2^{c_2} = \alpha^{(c_1,c_2)} \lesssim \sum_{I \in \I_0} \alpha^{\deg(I)}.$$
Thus to prove \eqref{targ} it will suffice to show
\be{targ-reduced}
|\Omega| \gtrapprox \alpha_2^{C/N}
\sum_{I \in \I_0} \alpha^{\deg(I)}.
\end{equation}

From \eqref{alpha-def} we expect, heuristically, that for generic points in $x \in \Omega$ and $t = O(1)$, the point $e^{tX_j}(x)$ has a probability $\sim \alpha_j$ of remaining in $\Omega$.  To make this precise we will need to pass from $\Omega$ to various refinements of $\Omega$, which we now construct.

We first make a technical reduction to reduce the diameter of $\Omega$.  By the pigeonhole principle (covering $V$ by balls of radius $O(\alpha_2^{2/N})$) one can find a subset ${\tilde \Omega} \subset \Omega$ which has diameter $O(\alpha_2^{2/N})$ and is such that
\be{om-0}
|{\tilde \Omega}| \gtrapprox \alpha_2^{C/N} |\Omega|.
\end{equation}
Henceforth we fix this ${\tilde \Omega}$.

\begin{definition}  Let $j$ be an integer.  A \emph{$j$-sheaf $\Omega'$ of width $w_j$} is defined to be a refinement $\Omega'$ of ${\tilde 
\Omega}$ such that, for all $x \in \Omega'$, the set
\be{times}
\{ |t| \ll 1: e^{tX_j}(x) \in \Omega' \}
\end{equation}
is a central set of width $w_j$ and measure $\gtrapprox \alpha_2^{C/N} \alpha_j$.
\end{definition}

This definition and the following lemma is very reminiscent of the restricted weak-type reductions used in, e.g. \cite{ccc}.  Our main innovation here is that the sets of times are central sets with fixed width; this will be crucial in placing (a large subset of) $\Omega$ inside one of the Carnot-Carath\'eodory balls used earlier.

\begin{lemma}\label{ref-lemma}  For any refinement $\Omega'$ of ${\tilde \Omega}$ and any integer $j$, there exists a refinement $\langle \Omega'\rangle_j$ of $\Omega'$ such that $\langle \Omega'\rangle_j$ is a $j$-sheaf of some width $\alpha_2^{C/N} \alpha_j \lessapprox w_j \lesssim \alpha_2^{2/N}$.
\end{lemma}

\begin{proof}
Fix $j$, $\Omega'$.  It will be convenient to change co-ordinates so that $\Sigma$ is locally like $\R^n$, and $X_j$ is just the constant vector field $e_n$.  In this case we can reparameterize $M_j$ as $\R^{n-1}$, so that $\pi_j: \R^n \to \R^{n-1}$ just becomes the Euclidean projection
$$ \pi_j(x_1,\ldots,x_n) = (x_1, \ldots, x_{n-1}).$$
We can of course let $0_\Sigma$ be the usual origin $(0, \ldots, 0)$.

Consider the function $(\pi_j)_* \chi_{\Omega'}$, the push-forward of $\chi_{\Omega'}$ by $\pi_j$.  This function has $L^1$ norm $|\Omega'|$, and is supported on a set of measure 
$$ \leq |\pi_j(\Omega)| = |\Omega|/\alpha_j \lessapprox \alpha_2^{-C/N} |\Omega'|/\alpha_j.$$
Because of this, we can find a set $E$ in $\R^{n-1}$ such that
$(\pi_j)_* \chi_{\Omega'} \gtrapprox \alpha_2^{C/N} \alpha_j$ on $E$
and
$$ \int_E (\pi_j)_* \chi_{\Omega'} \approx |\Omega'|.$$

At this point one could attempt to set $\langle \Omega' \rangle_j$ equal to $\Omega' \cap \pi_j^{-1}(E)$, but we would not get the crucial property that \eqref{times} is central with a fixed width.  To obtain this additional property we must refine further.

Let $x$ be any point in $E \subset \R^{n-1}$, and consider the set
$$ S(x) := \{ t: (x,t) \in \Omega' \}.$$
From the diameter restriction on $\Omega_0$ we see that this is a subset of $[-C\alpha_2^{2/N}, C\alpha_2^{2/N}]$ with measure $\gtrapprox \alpha_2^{C/N} \alpha_j$.

Let $I(x)$ be the dyadic interval of minimal length such that
$$ |I(x) \cap S(x)| \geq \frac{1}{4} |I(x)|^\eps |S(x)|;$$
if there are many intervals of this length, we choose $I(x)$ arbitrarily.  Let $w(x)$ denote the length of $I(x)$.

This interval is well-defined and must have length $\alpha_2^{C/N} \alpha_j \lessapprox w(x) \lesssim \alpha_2^{2/N}$.  Set $S'(x) := I(x) \cap S(x)$.  By construction we see that
\be{s-1}
|S'(x)| \approx |S(x)| \gtrapprox \alpha_2^{C/N} \alpha_j 
\end{equation}
and that
\be{s-2}
|I \cap S'(x)| \lesssim (|I|/w(x))^\eps |S'(x)|
\end{equation}
for all intervals $I$.  (Note that any interval $I$ can be covered by $O(1)$ dyadic intervals of comparable length, and each of the dyadic intervals are either contained in $I(x)$ or are disjoint from it.).

Since $w(x)$ can only take $O(\log(1/\alpha_2)) \approx 1$ many values, there exists an $\alpha_2^{C/N} \alpha_j \lessapprox w_j \lesssim \alpha_2^{2/N}$ such that 
$$ \int_{x \in E: w(x) = w_j} |S'(x)| \approx \int_E |S'(x)|,$$
which by the previous implies that
$$ \int_{x \in E: w(x) = w_j} |S'(x)| \approx \int_E |S(x)| =
\int_E (\pi_j)_* \chi_{\Omega'} \approx |\Omega'| \approx |{\tilde \Omega}|.$$
We then choose such a $w_j$ and define
$$ \langle \Omega' \rangle_j := \{ (x,t): x \in E, w(x) = w_j, t \in S'(x) \}.$$

We now verify that $\langle \Omega' \rangle_j$ is a $j$-sheaf of width $w_j$. 
Choose $(x_0, t_0) \in \langle \Omega' \rangle_j$.  By construction, $t_0$ lies in $S'(x_0)$, which is itself contained in an interval $I(x_0)$ of length $w_j$.  Since the set of times \eqref{times} corresponding to $(x_0,t_0)$ is just
$S'(x_0) - t_0$, we see from \eqref{s-1}, \eqref{s-2} that \eqref{times} is indeed a central set with width $w_j$ and measure $\gtrapprox \alpha_2^{C/N} \alpha_j$.  (The all-important centrality condition follows since $I(x_0)$ has length $w_j$, so that $S'(x_0) - t_0$ is contained in the interval $[-w_j, w_j]$.)
\end{proof}

The above Lemma allows one to refine a set $\Omega$ so that the $\pi_1$ fibers (for instance) have good properties.  It is tempting to iterate this lemma so that the $\pi_2$ fibers are also good, but unfortunately the latter refinement can destroy some of the properties of the former.  One could then try to iterate the lemma repeatedly, but the widths $w_1$, $w_2$ may change in doing so.  Fortunately we can salvage matters by using a pigeonholing argument (dating back at least to \cite{wolff:xray}; we use the formulation in \cite{laba:xray}.  Somewhat similar ideas appear in \cite{ccc}), if we are willing to lose a power of $\alpha_2^{1/N}$:

\begin{corollary}  There exists dyadic numbers $\alpha_2^{C/N} \alpha_j \lessapprox \mu_j \lesssim \alpha_2^{2/N}$ for $j=1,2$ and an increasing sequence
$$ \Omega_0 \prec \Omega_1 \ldots \prec \Omega_n \prec {\tilde \Omega}$$
such that for all $0 \leq j \leq n$, the set $\Omega_j$ is a $j$-sheaf whose width $w_j$ satisfies
$ \alpha_2^{C/N} \mu_j \lessapprox w_j \lesssim \mu_j$.  Here we have adopted the convention that the subscripts of $\mu_1$, $\mu_2$ are extended periodically.
\end{corollary}

\begin{proof}
We first define a much longer decreasing sequence of refinements
$$ {\tilde \Omega} = \Omega^{(0)} \succ \Omega^{(1)} \ldots \succ \Omega^{(N^2)}$$
by setting $\Omega^{(j)} = \langle \Omega^{(j-1)} \rangle_j$ for all $1 \leq j \leq N^2$.  (Note we are extending the construction in Lemma \ref{ref-lemma} periodically in $j$ in the usual fashion).  By construction we see that each $\Omega^{(j)}$ is an $j$-sheaf of width $w^{(j)}$ for some 
$\alpha_2^{C/N} \alpha_j \lessapprox w^{(j)} \lesssim \alpha_2^{2/N}$.  (Note that we are allowed to apply the lemma $N^2$ times because our implicit constants are allowed to depend on $N$).

Since $\Omega^{(j+2)} \subset \Omega^{(j)}$, we see from Corollary \ref{width-mono} that $w^{(j+2)} \lesssim w^{(j)}$ for all $1 \leq j \leq N^2-2$.  Thus if we set $N$ to be large enough, then by the pigeonhole principle one can find an even number $n \leq j_0 \leq N^2$ and numbers $\mu_1$, $\mu_2$ (extended periodically) such that 
$\alpha_2^{C/N} \alpha_j \lessapprox \mu_j \lesssim \alpha_2^{2/N}$ and $\alpha_2^{C/N} \mu_j \lessapprox w^{(j)} \lesssim \mu_j$ for all $j_0 - n \leq j \leq j_0$.  The claim then follows by setting $\Omega_j = \Omega^{(j_0 - j)}$ and $w_j = w^{(j_0 - j)}$ (note that $j_0 - j$ and $j$ always have the same parity).
\end{proof}

Let $\mu_j$, $\Omega_j$ be as above.  Since $\Omega_0$ is a refinement of ${\tilde \Omega}$, it is non-empty.  We fix a point $x_0 \in \Omega_0$.
Define $\delta_1, \delta_2$ by
$$ \delta_j := \alpha_2^{-1/N} \mu_j$$
and extend this periodically in the usual manner, $\delta_{j+2} = \delta_j$.
Then we have
$$ \alpha_2^{C/N} \delta_j \lessapprox w_j \lesssim \alpha_2^{1/N} \delta_j \hbox{ for all } j = 1, \ldots, n,$$
$$ 0 < \alpha_i/\delta_i \lessapprox \alpha_2^{-C/N},$$
and
$$ \alpha_2^{1+C/N} \lessapprox \delta_1, \delta_2 \lessapprox \alpha_2^{1/N}.$$
In particular we have
$$ \delta_1^{1/\eps} \lesssim \delta_2 \lesssim \delta_1^\eps.$$

For all $0 \leq i \leq n$ define the map $\Phi_i$ from a ball of radius $O(1)$ and center 0 in $\R^i$ to $\Sigma$ by
$$ \Phi_i(t_1, \ldots, t_i) := e^{t_i \delta_i X_i} \ldots e^{t_1 \delta_1 X_1} (x_0),$$
and let $T_i \subset \R^i$ denote the set
$$ T_i := \{ t \in \R^i: \Phi_j(t_1, \ldots, t_j) \in \Omega_{j-1} \hbox{ for all } 1 \leq j \leq i \},$$
with the convention that $T_0 = \R^0 = \{ () \}$, the singleton set consisting of the unique 0-tuple.  Observe that for each $0 \leq i < n$ and $t \in T_i$ the set
$$ \tau_i(t) := \{ t_{i+1} \in \R: (t, t_{i+1}) \in T_{i+1} \}$$
is a central set of measure $\gtrapprox \alpha_i/\delta_i$ and width $\tilde w_i$ for some
$\alpha_2^{C/N} \lessapprox \tilde w_i \lesssim \alpha_2^{1/N}$;
indeed, the sets $\tau_i(t)$ are just the sets \eqref{times} rescaled by $\delta_i$.

Since the $t_i$ are bounded (and this is the key place where we use the centrality of the set \eqref{times}), we see from \eqref{ywm-smooth} that 
$$ \Phi_j(T_j) \subset \Phi(B) \hbox{ for all } 1 \leq j \leq n$$
where $\Phi$, $B$ are as in Section \ref{balls-sec}; here we assume the radius of $B$ to be sufficiently large.

Now by construction, the set $\Phi_n(T_n)$ lies in $\Omega$.  By \eqref{volume} we thus have
$$ |\Omega| \geq |\Phi_n(T_n)| \approx |\Lambda(x_0)| |\Phi^{-1} \circ \Phi_n(T_n)|.$$
However, we see that
$$ \Phi^{-1} \circ \Phi_j = \Psi_j$$
where $\Psi_j$ maps a ball in $\R^j$ to $B$ and is defined by
$$ \Psi_j(t_1, \ldots, t_j) := e^{t_j Y_j} \ldots e^{t_1 Y_1} (0)$$
and the $Y_i$ are the vector fields from Section \ref{balls-sec}.  Thus we have
$$
|\Omega| \gtrapprox |\Lambda(x_0)| |\Psi_n(T_n)|.
$$
In the next section we shall show the bound
\be{psi-bound}
|\Psi_n(T_n)| \gtrapprox \alpha_2^{C/N} \frac{\alpha_1}{\delta_1} \ldots \frac{\alpha_n}{\delta_n}.
\end{equation}
Inserting this into the previous and using the fact that $\alpha_j$ and $\delta_j$ are periodic, we can write this as
$$
|\Omega| \gtrapprox \alpha_2^{C/N} |\Lambda(x_0)|
(\frac{\alpha_1}{\delta_1})^{\lfloor (n+1)/2 \rfloor}
(\frac{\alpha_2}{\delta_2})^{\lfloor n/2 \rfloor}.
$$
Now for all $I \in \I_0$, we observe that $\deg(I)_1, \deg(I)_2 \geq n-1$.  This is because all the vector fields $X_w$ in $I$ must be distinct (otherwise $\lambda_I \equiv 0$), and apart from the vector fields $X_1$, $X_2$, every other vector field $X_w$ has $\deg(w)_1, \deg(w)_2 \geq 1$.  Since\footnote{Notice the ``slack'' in the argument here.  What this is saying is that when one has a two-ends condition on the set of times \eqref{times}, then the lower bounds on $|\Omega|$ improve substantially.  This is consistent with the experience with the Kakeya problem in e.g. \cite{wolff:kakeya}.} 
$$n-1 \geq \lfloor (n+1)/2 \rfloor, \lfloor n/2 \rfloor$$
and $\alpha_j/\delta_j \lessapprox 1$, we thus have
$$
|\Omega| \gtrapprox \alpha_2^{C/N} |\Lambda(x_0)|
(\frac{\alpha_1}{\delta_1})^{\deg(I)_1}
(\frac{\alpha_2}{\delta_2})^{\deg(I)_2}
$$
for all $I \in \I_0$.  Combining this with \eqref{lambda-bound} we obtain \eqref{targ-reduced} as desired.

It remains to show \eqref{psi-bound}.  This will occupy the entirety of the next section.

\section{Conclusion of the argument}\label{conclude-sec}

We now prove \eqref{psi-bound}.  To simplify the notation we write $\beta_j := \alpha_j / \delta_j$.  Our situation is as follows.  We have a point $x_0 \in V$ which is within $O(1/K)$ of the origin 0.  We have sets $T_i \in \R^i$ for $i = 0, \ldots, n$, which have the structure
$$ T_{i+1} = \{ (t, t_{i+1}) \in \R^{i+1}: t \in T_i; t_{i+1} \in \tau_i(t) \}$$
for $i = 0, \ldots, n-1$, with $T_0 = \R^0 = \{ () \}$ and for each $t \in T_i$, the set $\tau_i(t)$ is a central set of width 
$$ \alpha_2^{C/N} \lessapprox {\tilde w}_i \lesssim \alpha_2^{1/N}$$
and $\tau_i(t)$ has measure 
\be{beta-bound}
|\tau_i(t)| \gtrapprox \alpha_2^{C/N} \beta_i.
\end{equation}
From our upper bound on $\tilde w_i$ we see that
$$ T_i \subseteq [-r,r]^i$$
for all $0 \leq i \leq n$ and some $r \sim \alpha_2^{1/N}$, which we fix.

From \eqref{beta-bound} and Fubini's theorem we have
$$ |T_{i+1}| \gtrapprox \alpha_2^{C/N} \beta_i |T_i|$$
and hence
$$ |T_n| \gtrapprox \alpha_2^{C/N} \beta_1 \ldots \beta_n.$$
On the other hand, we wish to prove that
\be{reduced}
|\Psi_n(T_n)| \gtrapprox \alpha_2^{C/N} \beta_1 \ldots \beta_n.
\end{equation}
It will therefore suffice to find a subset $T'_n$ of $T_n$ with
\be{large-tpn}
|T'_n| \gtrapprox \alpha_2^{C/N} \beta_1 \ldots \beta_n
\end{equation}
such that we have the Jacobian bound
\be{implicit}
|\det(D\Psi_n)(t)| \gtrapprox \alpha_2^{C/N} \hbox{ for all } t \in T'_n.
\end{equation}
Indeed, if we can find such a $T'_n$, then we cover $T'_n$ by balls of radius $\alpha_2^{C_1/N}$ for some large constant $C_1$.  If $C_1$ is large enough, then
$\Psi_n$ is one-to-one on every one of the balls that intersects $T'_n$, thanks
to \eqref{implicit}, the smoothness of $\Psi_n$ (recall from Section \ref{balls-sec} that the $\Psi_n$ are smooth uniformly in $x_0$, $\delta_1$, $\delta_2$), and the inverse function theorem.  By the pigeonhole principle one of these balls must contain at least $\alpha_2^{CC_1/N}$ of the set $T'_n$, and the claim \eqref{reduced} then follows from \eqref{implicit}, \eqref{large-tpn} and the change-of-variables formula.
                              
It remains to construct a set \eqref{large-tpn}.  This will be achieved by the width property of the $\tau_i(t)$ and several applications of Lemma \ref{lower},
but first we need to control some derivative of $\det(D\Psi_n(t))$.  We shall do so using the H\"ormander condition assumption, by invoking the following quantitative analogue of Frobenius's theorem (which can also be derived from the arguments in \cite{cnsw}):

\begin{lemma}\label{mixed}  There exists a multi-index $\beta = (\beta_1, \ldots, \beta_n)$ with $|\beta| \leq C(I_0,n)$ such that
$$ |\partial_t^\beta \det(D\Psi_n)|_{t=0}| \gtrsim C_K^{-1} .$$
\end{lemma}

\begin{proof}
We shall prove a more general statement.  For any $1 \leq i \leq n$, let $J_i$ denote the $i$-form
$$ J_i(t) := \partial_{t_1} \Psi_i(t) \wedge \ldots \wedge \partial_{t_i} \Psi_i(t).$$
defined for all $t$ in a ball in $\R^i$.  We shall show that for each $i$ there exists a multi-index $\beta^i$ with $|\beta^i| \leq C(I_0,n)$ and a small number $0 < \theta_i \ll 1$ depending on $K$ such that
\be{j-bound}
|\partial_t^{\beta^i} J_i|_{t=0}| \gtrsim \theta_i.
\end{equation}
The Lemma will then follow by applying \eqref{j-bound} with $i=n$.

The idea of the proof can be sketched as follows.
Informally, the statement \eqref{j-bound} asserts that the vector fields $X_1, \ldots, X_i$ are non-degenerate, and so the iterated map $\Psi_i$ fills out a large subset of a $j$-dimensional manifold.  If $X_{i+1}$ is transverse to this manifold then one can induct and obtain \eqref{j-bound} for $i+1$.  Similarly if $X_{i+1}$ is tangent to the manifold to some (bounded) finite order.  The only remaining case is if $X_{i+1}$ is tangent to the manifold to some extremely high order.  But this will imply that the iterated Lie brackets of $X_{i+1}$ and $X_i$ (for instance) are also tangent to the manifold, and this will eventually contradict the H\"ormander condition assumption.

We now return to the rigorous proof of \eqref{j-bound}.  For $i=1$ this is clear just by taking $\beta = 0$ and using the non-vanishing of $Y_1$.  Now suppose that $i>1$, and the claim has already been proven for $i-1$.

Let $0 < \theta_i \ll 1$ be a small number depending on $K$, $\theta_{i-1}$ to be chosen later, let $A = A(I_0,n)$ be a large integer to be chosen later, and let $\epsilon = \epsilon(A)$ be a small number to be chosen later.

By induction hypothesis and the smoothness of $J_{i-1}$, we can find a ball $B_{i-1} \subset \R^{i-1}$ of radius $\gtrsim C_{\epsilon,\theta_{i-1}}^{-1}$ and at a distance $O(\epsilon \theta_{i-1})$ from the origin such that 
\be{j-bo}
|J_{i-1}(t)| \gtrsim \epsilon^C \theta_{i-1}^C \hbox{ for all } t \in B_{i-1}
\end{equation}
(cf. the arguments in Lemma \ref{lower}); the constant $C$ can depend on $\beta_{i-1}$.  By shrinking $B_{i-1}$ if necessary we may then assume that $\Psi_{i-1}$ is one-to-one on $B_{i-1}$.  The set $\Psi_{i-1}(B_{i-1})$ is thus an open subset of a smooth $i-1$-dimensional manifold in $B$.  Observe that the $i-1$-form $J_{i-1}(t)$ is tangent to $\Psi_{i-1}(B_{i-1})$ at $\Psi_{i-1}(t)$ for every $t \in B_{i-1}$.

Suppose for contradiction that we have
$$ |\partial_t^\beta J_i|_{t=0}| \lesssim \theta_i$$
for all multi-indices $\beta$ with $|\beta| \leq A$.  Then by Taylor expansion we have
$$ |J_i(t_1, \ldots, t_{i-1},0)| \lesssim C_A \theta_i + C_A \epsilon^A \theta_{i-1}^{A}$$
whenever $t_1, \ldots, t_{i-1} = O(\epsilon \theta_{i-1})$, and more generally that
$$ |\partial_t^\beta J_i(t_1, \ldots, t_{i-1},0)| \lesssim C_A \theta_i + C_A (\epsilon\theta_{i-1})^{A-|\beta|}$$
whenever $|\beta| < A$.  On the other hand, by definition we have
\bas
J_i(t_1, \ldots, t_{i-1},0) 
&= J_{i-1}(t_1,\ldots,t_{i-1}) 
\wedge \partial_{t_i} \Psi_i(t_1,\ldots,t_{i-1},t_i)|_{t_i=0}\\
&= J_{i-1}(t_1,\ldots,t_{i-1}) 
\wedge Y_i(\Psi_{i-1}(t_1,\ldots,t_{i-1})).
\end{align*}
Because of this and \eqref{j-bo}, we see that for every $t \in B_{i-1}$ that we can write
$$ Y_i(\Psi_{i-1}(t)) = Y^{tangent}(\Psi_{i-1}(t)) + O(C_{A,\theta_{i-1},\epsilon} \theta_i + C_A (\epsilon\theta_{i-1})^{A-C})$$
where $Y^{tangent}(\Psi_{i-1}(t))$ is tangent to $\Psi_{i-1}(B_{i-1})$ at $\Psi_{i-1}(t)$.  Since $\Psi_{i-1}(B_{i-1})$ is a smooth manifold, we may thus write
$$ Y_i = Y^{tangent} + Y^{error}$$
on a $C_{\theta_{i-1}}^{-1}$-neighborhood of $\Psi_{i-1}(B_{i-1})$, where $Y^{tangent}$ is a smooth vector field which when restricted to $\Psi_{i-1}(B_{i-1})$ is tangent to $\Psi_{i-1}(B_{i-1})$, and $Y^{error}$ is another smooth vector field such that
$$ |\nabla^\beta Y^{error}| \lesssim C_{A,\theta_{i-1},\epsilon} \theta_i + C_A (\epsilon\theta_{i-1})^{A-|\beta|-C}$$
for all multi-indices $\beta = (\beta_1, \ldots, \beta_n)$ with $|\beta| \ll A$.

From the identity
$$ Y_{i-1}(\Psi_{i-1}(t)) = \partial_{t_{i-1}} \Psi_{i-1}(t)$$
we see that $Y_{i-1}$ is also tangent to $\Psi_{i-1}(B_{i-1})$ when restricted to $\Psi_{i-1}(B_{i-1})$.  Thus any Lie bracket combination of $Y_{i-1}$ and $Y^{tangent}$ must also be tangent to $\Psi_{i-1}(B_{i-1})$ when restricted to $\Psi_{i-1}(B_{i-1})$.

Now let $w$ be any word with $\deg(w)_1 + \deg(w)_2 \ll A$.  The vector field $Y_w$ is $K^{1- \deg(w)_1 - \deg(w)_2}$ times a Lie bracket combination of $Y_1$ and $Y_2$.  One of these vector fields is $Y_{i-1}$, the other is $Y_i$, which can be split into $Y^{tangent}$ and $Y^{error}$.  From the above discussion we thus see that $Y_w$ when restricted to $\Psi_{i-1}(B_{i-1})$ is equal to a vector field tangent to $\Psi_{i-1}(B_{i-1})$, plus an error term of magnitude at most
$$ C_{K,A,\theta_{i-1},\epsilon} \theta_i + C_{K,A} \epsilon \theta_{i-1}.$$
Applying this with $Y_{w_1}, \ldots Y_{w_n}$ and taking wedge products we obtain
$$ |Y_{w_1}(\Psi_{i-1}(t)) \wedge \ldots \wedge Y_{w_n}(\Psi_{i-1}(t))| \lesssim C_{K,A,\theta_{i-1},\epsilon} \theta_i + C_{K,A} \epsilon \theta_{i-1}.$$
If we pick $\epsilon$ sufficiently small depending on $A$ and $K$ and $\theta_i$ sufficiently small depending on $A, K$, $\theta_{i-1}$, $\eps$ we see that this bound contradicts \eqref{det-large}, and we are done.
\end{proof}

We can now quickly construct a set $T'_n$ with the properties \eqref{large-tpn}, \eqref{implicit}, which will finish the proof of \eqref{psi-bound} and thus of
Theorem \ref{main}.

By a Taylor expansion we may write
$$ \det(D\Psi_n)(t) = P(t) + O(\alpha_2^N)$$
on $t \in [-r,r]^n$, where $P$ is a polynomial of degree $O(N^2) = O(1)$,
$\|P\|_\infty \lesssim 1$ and $|\partial_t^{\alpha} P (0)| \gtrsim
C_{K}^{-1} $ (the index $\alpha = (\alpha_1,\ldots, \alpha_n )$ with
bounded coefficients is as in Lemma \ref{mixed}). It therefore suffices
to construct $T_n^{'}$ satisfying \eqref{large-tpn} so that 
\eqref{implicit} holds with $P(t)$ in place of ${\rm det}(D\Psi_n)(t)$.

We apply Lemma \ref{lower} with $S$ equal to $\tau_1(())$.  This gives us a subset $\tau'_1(())$ of $\tau_1()$ with measure
$$ |\tau'_1(())| \gtrapprox \alpha_2^{C/N} \beta_1$$
such that
$$ | \partial_{t_2}^{\alpha_2} \ldots \partial_{t_n}^{\alpha_n} P(t_1,0,\ldots,0) | \gtrapprox \alpha_2^{C/N}$$
for all $t_1 \in \tau'_1(())$.
 
Fix $t_1 \in \tau'_1(())$.  We now apply Lemma \ref{lower} with $S$ equal to $\tau_2(t_1)$.  This gives us a subset $\tau'_2(t_1)$ of $\tau_2(t_1)$ with measure$$ |\tau'_2(t_1)| \gtrapprox \alpha_2^{C/N} \beta_2$$
such that
$$ | \partial_{t_3}^{\alpha_3} \ldots \partial_{t_n}^{\alpha_n} P(t_1,t_2,0,\ldots,0) | \gtrapprox \alpha_2^{C/N}$$
for all $t_2 \in \tau'_2(t_1)$.
 
We continue in this manner, defining sets $\tau'_i(t_{i-1})$ recursively.  Note
that since our constants $C$ are allowed to depend on $n$ that we will not have
trouble iterating the above scheme $n$ times.  If we recursively define $T'_0 =
\{()\}$ and
$$ T'_{i+1} := \{ (t, t_{i+1}) \in \R^{i+1}: t \in T'_i; t_{i+1} \in \tau'_i(t)
\}$$                          
we thus see that $T'_n$ obeys \eqref{large-tpn} and \eqref{implicit} as desired.  This completes the proof of Theorem \ref{main}.

\section{Remarks}\label{remarks-sec}

\begin{itemize}

\item The estimate \eqref{targ-reduced} can be rewritten using \eqref{heur}, and the fact that $N$ was arbitrary, in the rather appealing form
\be{iso}
|\Omega| \gtrapprox |B(0; \alpha_1, \alpha_2)|.
\end{equation}
In other words, if one fixes the thicknesses $\alpha_j = |\Omega|/|\pi_j(\Omega)|$ of a set $\Omega$ in a small neighbourhood of 0, then the volume of the set $\Omega$ is essentially minimized when $\Omega$ is a two-parameter ball centered at the origin.  One may thus think of our main result as a kind of two-parameter isoperimetric inequality in the spirit of \cite{loomis:whitney}.  

\item It would be nice if the $\gtrapprox$ in \eqref{iso} could be replaced by a $\gtrsim$, as this would give restricted-type boundedness on the entire boundary of the Newton polytope $P$ (bringing these results in line with those in, e.g. \cite{ccc}).  If one wishes to adapt the above argument to do this, there seem to be two major obstacles.  The first is that the argument requires more and more regularity on $(\Sigma, M_1, M_2, \pi_1, \pi_2, 0_\Sigma)$ as one approaches the boundary of $P$.  This particular difficulty can be avoided by restricting one's attention to model cases, such as when $\Sigma$ is a nilpotent Lie group and $X_1$, $X_2$ are left-invariant vector fields.  The other major difficulty is that one cannot refine the set of times \eqref{times} to a central set of a fixed width without losing an $\eps$ in the exponents.  We do not know how to get around this loss.

\item One can be more ambitious still, and try to obtain strong-type boundedness at the endpoints of $P$.  This is already difficult to do (at least if one only uses geometric combinatorics as in this paper) in such model cases as convolution with a compactly supported measure on the parabola $\{(t, t^2)\}$ in $\R^2$.  In this model case the endpoint mapping is $L^{3/2} \to L^3$, but a naive combinatorial argument only gives $L^{3/2,1} \to L^{3,\infty}$ (however see
\cite{ober:circle}; also see \cite{ccc} for more discussion on this).  Of course, one can obtain the strong-type endpoint in this case if one is willing to use such tools as complex interpolation and $L^2$ smoothing estimates (see e.g. \cite{littman}).  It may also be possible to do so by pure geometric combinatorics, but one probably has to control the extent to which various two-parameter balls of varying radii can overlap each other.  In particular, a two-parameter covering lemma of some sort may be needed.

\item It seems likely that one could extend these results (with the aid of the results in \cite{cnsw} or \cite{greenblatt}) to obtain good $(L^{p_1}, L^{p'_2})$ mapping properties for 
fractional integral operators such as
$$ T_\alpha f(x) := \int f(\gamma(x,t)) \frac{dt}{|t|^\alpha}$$
for $0 < \alpha < 1$, where for each $x$, the map $t \mapsto \gamma(x,t)$ is a smooth parameterization of the curve 
$\gamma_x$ of the type discussed in the introduction, and the H\"ormander condition obeyed. In fact, if $\gamma(x,0) \equiv x$ and
$\alpha = 1 - \frac{1}{k}$ with $k$ odd, then one easily sees (from a simple change of variables) that $T_{\alpha}$ falls within 
the class of operators described by \eqref{duality} where the corresponding vector fields $X_1$ and $X_2$ are now linearly
dependent at the origin. 

\item The methods here should probably be able to give an alternate way to obtain the $L^p$ boundedness of maximal operators of the form 
$$ M f(x) := \sup_{k \geq 0} 2^{k} \int f(\gamma(x,t)) a(2^k t)\ dt;$$
in particular, it is now quite possible that one could prove boundedness of such operators using purely geometric methods (avoiding the Fourier transform, as is used in e.g. \cite{cnsw}).
This is of course not the only type of maximal function one could study; more generally one could consider operators of the form
$$ M f(x) := \sup_{r} \int f(\gamma(x,t,r)) a(t)\ dt$$
for some parameter $r$; model examples here include the circular maximal function and the Kakeya maximal function (see \cite{wolff:survey} for a survey of these operators).  It seems difficult however to adapt our techniques to these operators except in very low dimensions; one would have to linearize the parameter $r$ as $r(x)$, and this begins to destroy the differentiability properties of the iterated flow map $\Phi_n$ when $n > 3$.  On the other hand, Christ and Erdogan have recently used geometric combinatorics ideas of a similar flavor to those in \cite{ccc} and in this paper to obtain sharp mixed-norm estimates for certain classes of x-ray transforms; see \cite{ce}.  However, the
general question of understanding maximal operators of this type is still far beyond our current technology.

\item It is also tempting to try to use these techniques to prove $L^p$ smoothing estimates (i.e. mapping $L^p$ to $L^p_\alpha$ rather than to $L^q$).  However this problem is much more difficult due to the presence of cancellation, and does not have a nice geometric interpretation such as \eqref{iso}.  Indeed this problem is quite difficult even for model cases such as convolution with the curve $(t,t^2,t^3)$, as it is related to the local smoothing conjecture (see \cite{smith:helix}).  On the other hand, Wolff \cite{wolff:smsub} has recently combined geometric combinatorics techniques with Fourier methods to obtain some progress on these types of problems. However, with respect to curves in the plane, Seeger \cite{seeger:jams}, \cite{seeger:duke} has obtained
sharp (up to endpoints) $L^p$ to $L^p_\alpha$ estimates. 

\item  Another possible generalization would be to higher-dimensional averages, or to asymmetric averages in which the manifolds $M_1$ and $M_2$ have different dimensions.  In our language, this would mean that $\pi_1$ and $\pi_2$ now have corank $k_1$ and $k_2$ which are possibly greater than 1, and we would replace the single vector field $X_j$ by a family of commuting vector fields $X_j^1, \ldots, X_j^{k_j}$.  While it is possible to use this machinery (perhaps combined with the techniques in \cite{cnsw}) to get some non-trivial $(L^p,L^q)$ mapping result, it seems difficult to obtain sharp results, because there seems to be no satisfactory analogue of the notion of width in more than one dimension. See Seeger \cite{seeger:jams} where nontrivial $(L^p,L^q)$ as well as
$(L^p,L^p_\alpha)$ results are obtained when $\Sigma$ is a hypersurface in $M_1 \times M_2$ (i.e., dim$\Sigma = k_1 + k_2 + 1$).  

\item It is crucial in our arguments that the sets $\tau_i(t)$ of times are central, so that the times $t_i$ are always close to 0.  This allows us to restrict the set $\phi_n(T_n)$ to lie inside a small two-parameter Carnot-Carath\'eodory ball, whose geometry can be well controlled by the machinery of Section \ref{balls-sec}.  If we allowed $\tau_i(t)$ to wander far away from the origin, then it would in fact be impossible to obtain the estimate
\be{phin}
|\Phi_n(T_n)| \gtrapprox \alpha_1^{c_1} \alpha_2^{c_2}
\end{equation}
for the desired values of $c_1, c_2$, even in the model case when the $\tau_i(t)$ are all intervals of length $\alpha_i$ and we are considering convolution with $(t, \ldots, t^{n-1})$ for some $n > 5$ (see \cite{ccc} for some further discussion of this issue; this was also independently observed by Greenblatt).  The point is that we lose control of the geometry if one flows too far along one or more vector fields.  On the other hand, failure of the estimate \eqref{phin} does not imply that failure of the lower bound on $|\Omega|$, because $\Phi_n(T_n)$ may only occupy a small portion of $\Omega$.  Our particular selection method for $T_n$ (using the machinery of $j$-sheaves) is thus essential to ensure that $\Phi_n(T_n)$ does not degenerate to only a small fraction of $\Omega$.

\end{itemize}

\end{document}